\documentclass[11pt]{article}
\usepackage{epsfig}
\usepackage{multirow}
\usepackage{color}
\usepackage{amsmath}
\usepackage{amssymb}
\makeatletter
\def\singlespace{\def\baselinestretch{1}\@normalsize}

\@addtoreset{equation}{section}
\renewcommand{\baselinestretch}{1.412}
\renewcommand{\theequation}{\arabic{section}.\arabic{equation}}

\newcommand{\convD}{\stackrel{D}{\longrightarrow}}

\usepackage{natbib}
\allowdisplaybreaks

\newcommand{\T}{\!\mbox{\scriptsize T}}
\newcommand{\bX}{\mbox{\bf X}}
\newcommand{\0}{\mbox{\bf 0}}
\newcommand{\CV}{\mbox{CV}}

\newcommand{\MSE}{\mbox{MSE}}
\newcommand{\MISE}{\mbox{MISE}}
\newtheorem{thm}{Theorem}

\newtheorem{rmk}{Remark}

\def\wt{\widetilde}

\def\wh{\widehat}
\def\wb{\overline}

\def\log{\hbox{log}}

\def\boxit#1{\vbox{\hrule\hbox{\vrule\kern6pt
          \vbox{\kern6pt#1\kern6pt}\kern6pt\vrule}\hrule}}

\def\var{\hbox{var}}

\def\bse{\begin{eqnarray*}}
\def\ese{\end{eqnarray*}}
\def\be{\begin{eqnarray}}
\def\ee{\end{eqnarray}}
\def\bq{\begin{equation}}
\def\eq{\end{equation}}

\def\wh{\widehat}

\def\T{{\cal T}}

\def\trans{^{\rm T}}

\newtheorem{proposition}{Proposition}

\newcommand{\blem}{\begin{lemma}}
\newcommand{\elem}{\end{lemma}}
\newcommand{\bthe}{\begin{theorem}}
\newcommand{\ethe}{\end{theorem}}

\newtheorem{definition}{Definition}
\newtheorem{lemma}[definition]{Lemma}
\newtheorem{theorem}[definition]{Theorem}

\def\delete#1{\iffalse #1 \fi}

\def\bse{\begin{eqnarray*}}
\def\ese{\end{eqnarray*}}
\def\bee{\begin{enumerate}}
\def\eee{\end{enumerate}}
\def\bqe{\begin{eqnarray}}
\def\eqe{\end{eqnarray}}
\def\bed{\begin{description}}
\def\eed{\end{description}}
\def\bei{\begin{itemize}}
\def\eei{\end{itemize}}

\def\pmb#1{\setbox0=\hbox{#1}%
    \kern-.025em\copy0\kern-\wd0
    \kern.05em\copy0\kern-\wd0
    \kern-.025em\raise.0433em\box0 }
\def\pmbh#1#2{\setbox0=\hbox{#1}%
    \setbox1=\hbox{#2}%
    \kern-.025em\copy0\kern-\wd0
    \kern.05em\copy1\kern-\wd0
    \kern-.025em\raise.0433em\box0 }

\def\frac#1#2{{#1\over#2}}

\def\boxit#1{\vbox{\hrule\hbox{\vrule\kern6pt
   \vbox{\kern6pt#1\kern6pt}\kern6pt\vrule}\hrule}}

\def\listing#1{\vskip 4mm\begin{verbatim}\input#1 \vskip 4mm}
\def\thick#1{\hbox{\rlap{$#1$}\kern0.25pt\rlap{$#1$}\kern0.25pt$#1$}}

\def\wt{\widetilde}
\def\wh{\widehat}

\def\var{\mbox{var}}

\def\T{{\mbox{\rm\tiny T}}}

\def\pmbh{{\pmb h}}

\def\calM{{\cal M}}
\def\calN{{\cal N}}

\renewcommand\today{\ifcase\month\or
   Jan\or Feb\or Mar\or Apr\or May\or
   Jun\or Jul\or Aug\or Sep\or Oct\or Nov\or
   Dec\fi
   \space\number\day, \number\year}

\newcommand{\va}{{\bf a}}
\newcommand{\vb}{{\bf b}}

\newcommand{\ve}{{\bf e}}

\newcommand{\vE}{{\bf E}}
\newcommand{\vG}{{\bf G}}

\newcommand{\vK}{{\bf K}}

\newcommand{\vu}{{\bf u}}
\newcommand{\vv}{{\bf v}}

\newcommand{\vA}{{\bf A}}
\newcommand{\vB}{{\bf B}}
\newcommand{\vC}{{\bf C}}

\newcommand{\vI}{{\bf I}}
\newcommand{\vJ}{{\bf J}}

\newcommand{\vP}{{\bf P}}
\newcommand{\vQ}{{\bf Q}}

\newcommand{\vR}{{\bf R}}

\newcommand{\vU}{{\bf U}}
\newcommand{\vV}{{\bf V}}
\newcommand{\vX}{{\bf X}}

\newcommand{\vD}{{\bf D}}
\newcommand{\vO}{{\bf O}}

\newcommand{\vnull}{{\bf 0}}

\def\vPhi{\mathbf{\Phi}}
\def\vPsi{\mathbf{\Psi}}

\newcommand{\vdelta}{\mbox{\boldmath $\delta$}}

\newcommand{\vSigma}{\mbox{\boldmath $\Sigma$}}

\newcommand{\vTheta}{\mbox{\boldmath $\Theta$}}
\newcommand{\valpha}{\mbox{\boldmath $\alpha$}}
\newcommand{\vbeta}{\mbox{\boldmath $\beta$}}
\newcommand{\vgamma}{\mbox{\boldmath $\gamma$}}

\newcommand{\vPi}{\mbox{\boldmath $\Pi$}}

\newcommand{\veta}{\mbox{\boldmath $\eta$}}

\newcommand{\vmu}{\mbox{\boldmath $\mu$}}

\newcommand{\vtheta}{\mbox{\boldmath $\theta$}}
\newcommand{\vxi}{\mbox{\boldmath $\xi$}}

\newcommand{\bay}{\begin{array}}
\newcommand{\eay}{\end{array}}

\newcommand{\bqa}{\begin{eqnarray*}}
\newcommand{\eqa}{\end{eqnarray*}}
\newcommand{\bqan}{\begin{eqnarray}}
\newcommand{\eqan}{\end{eqnarray}}
\newcommand{\bqt}{\begin{quote}}
\newcommand{\eqt}{\end{quote}}
\newcommand{\bt}{\begin{tabbing}}
\newcommand{\et}{\end{tabbing}}
\newcommand{\bit}{\begin{itemize}}
\newcommand{\eit}{\end{itemize}}
\newcommand{\ben}{\begin{enumerate}}
\newcommand{\een}{\end{enumerate}}
\newcommand{\beq}{\begin{equation}}
\newcommand{\eeq}{\end{equation}}
\newcommand{\bdefi}{\begin{definition}}
\newcommand{\edefi}{\end{definition}}
\newcommand{\bpro}{\begin{proposition}}
\newcommand{\epro}{\end{proposition}}
\newcommand{\bco}{\begin{corollary}}
\newcommand{\eco}{\end{corollary}}
\newcommand{\bdes}{\begin{description}}
\newcommand{\edes}{\end{description}}

\newcommand{\vepsilon}{\mbox{\boldmath $\epsilon$}}

\title{\bf 
Homogeneity Pursuit in Single Index Models based Panel Data Analysis
}

\author{
Heng Lian
\\
Department of Mathematics
\\
City University of Hong Kong, Kowloon, Hong Kong
\and
Xinghao Qiao
\\
Department of Statistics
\\ 
London School of Economics, United Kingdom
\and
Wenyang Zhang
\\
Department of Mathematics
\\
The University of York, United Kingdom
}
\date{}

\begin{document}

\maketitle

\abstract{
Panel data analysis is an important topic in statistics and 
econometrics.  Traditionally, in panel data analysis, all individuals are 
assumed to share the same unknown parameters, e.g. the same coefficients 
of covariates when the linear models are used, and the differences between the 
individuals are accounted for by cluster effects.  This kind of modelling only 
makes sense if our main interest is on the global trend, this is because it 
would not be able to tell us anything about the individual 
attributes which are sometimes very important.  In this paper, we proposed 
a modelling based on the single index models embedded with homogeneity for 
panel data analysis, which builds the individual attributes in the model and is 
parsimonious at the same time.  We develop a data driven approach to identify 
the structure of homogeneity, and estimate the unknown parameters and 
functions based on the identified structure.  Asymptotic properties of the 
resulting estimators are established.  Intensive simulation studies conducted 
in this paper also show the resulting estimators work very well when sample 
size is finite.  Finally, the proposed modelling is applied to a public 
financial dataset and a UK climate dataset, the results reveal some interesting 
findings.
}

\vskip0.7cm






\section{Introduction}

\subsection{Preamble}
\label{pre}

Panel data analysis is an important topic in statistics and econometrics.  The 
traditional approach for analysing panel data assumes all individuals share the 
same unknown parameters, and uses cluster effects to account for the difference 
between individuals.  For example, when the linear models are 
used, the coefficients of the covariates are assumed to be the same across all 
individuals, i.e.
$$
y_{it} = \vX_{it}^{\T} \vbeta + \epsilon_{it},
\quad
i=1, \ \cdots, \ m;  \  \  t=1, \ \cdots, \ T,
$$
where $y_{it}$ and $\vX_{it}$, a $(p+1)$-dimensional vector, are respectively 
the $t$th observations of the response variable and covariate of the $i$th 
individual.  $\epsilon_{it}$, $t=1, \ \ldots, \ T$, are correlated for any 
given $i$, and the cluster effects are included in $\epsilon_{it}$.  See 
\cite{hsiao14} and the reference therein.  Whilst this modelling idea is 
useful when the global trend of the impact of a covariate on the response 
variable is of our main interest, it does not tell us anything about the 
individual attributes which are sometimes very important.

In order to explore the individual attributes, we need to make them more 
concrete and distinctive in modelling.  A simple approach to do so would be 
using
\begin{equation}
y_{it} = \vX_{it}^{\T} \vbeta_i + \epsilon_{it},
\quad
i=1, \ \cdots, \ m;  \  \  t=1, \ \cdots, \ T,
\label{eq1}
\end{equation}
to fit the data.  However, this modelling approach would result in $m(p+1)$ 
unknown coefficients to estimate, which is too many, because $m$ is 
usually of the magnitude of hundreds, or even more, in practice.  This 
modelling also ignores the similarity which may exist among some individuals.  
Such similarity may have very important practical meaning, and could lead to 
some important findings in practice.  In addition to that, statistically 
speaking, the modelling, like (\ref{eq1}) without any conditions imposed, would 
also pay a price on variance side of the estimators resulted because the 
available information is not used up.

In order to explore the individual attributes and account for the similarity 
among some individuals at the same time, \cite{ke15} proposed a penalised 
likelihood/least squares based approach to pursue the homogeneity in the 
linear models, i.e. (\ref{eq1}), used for panel data analysis, under the 
framework of treating homogeneity as a kind of sparsity.   Regression under 
homogeneity condition has also been studied by quite a few recent works, e.g. 
 \cite{tibshirani05,friedman07,bondell2008simultaneous,jiangqian13}
, and the references therein.  Like \cite{ke15}, 
the methods in these works are all based on penalised likelihood/least 
squares.  \cite{ke16} took a different approach, they formulated the 
homogeneity pursuit problem as a problem of change point detection and applied 
the binary segmentation approach to identify the homogeneity in the
linear models with interactive effects.

The existing literature about homogeneity pursuit mainly focuses on the linear 
models.  It is well known that the linearity condition may not hold for many 
datasets, and the exploration of linear relationship is not sufficient in many 
cases.  As a consequence, the semiparametric modelling is becoming more and 
more useful in panel data analysis.  Among various semiparametric models, the 
single index models 
have many advantages, and are a very successful tool in data analysis, see 
 \cite{hardle1989investigating,carroll97,yuruppert02,zhuxue06,xia08,penghuang11,zhumiaopeng12,guo2016dynamic}, 
 and the 
reference therein.  In this paper, we are going to investigate the homogeneity 
pursuit in the single index models used for panel data analysis.  The detailed 
definition of the models we are going to address in this paper is given in 
Section \ref{mo0}

\subsection{The single index models with homogeneity structure}
\label{mo0}

Let $y_{it}$ and $\vX_{it}$, a $(p+1)$-dimensional vector, be 
respectively the $t$th observations of the response variable and covariate of 
the $i$th individual, $i=1, \ \cdots, \ m; \ \  t=1, \ \cdots, \ T$. We 
consider the models 
\begin{equation}
y_{it}
=
g_i(\vX_{it}\trans\vbeta_i) + \epsilon_{it},
\quad 
i=1, \ \cdots, \ m; \ \  t=1, \ \cdots, \ T,
\label{model}
\end{equation}
where 
\begin{equation}
g_i(\cdot)
=
\left\{
\begin{array}{ll}
\ g_{(1)}(\cdot) & \mbox{ when }  i \in G_{1,1},
\\
\ g_{(2)}(\cdot) & \mbox{ when }  i \in G_{1,2},
\\
\quad \vdots & \quad \vdots
\\
\ g_{(H_1)}(\cdot) & \mbox{ when }  i \in G_{1,H_1},
\end{array}
\right.
\quad
\beta_{ij}
=
\left\{
\begin{array}{ll}
\ \beta_{(1)} & \mbox{ when }  (i, j) \in G_{2,1},
\\
\ \beta_{(2)} & \mbox{ when }  (i, j) \in G_{2,2},
\\
\quad \vdots & \quad \vdots
\\
\ \beta_{(H_2)} & \mbox{ when }  (i, j) \in G_{2,H_2},
\end{array}
\right.
\label{con}
\end{equation}
$\mathbb{G}_1=\{G_{1,k}: \ k = 1, \ \cdots, \ H_1\}$ is a partition of set 
$\{1, \ \cdots, \ m\}$, $\mathbb{G}_2=\{G_{2,k}: \ k = 1, \ \cdots, \ H_2\}$ is a partition 
of set $\{(i,j): \ i = 1, \ \cdots, \ m; \ j=1, \ \cdots, \ p\}$, $\beta_{ij}$ 
is the $(j+1)$th component of $\vbeta_i$, 
and
$$
E(\epsilon_{it}| \vX_{it}) = 0,
\quad
\var(\epsilon_{it}| \vX_{it}) = \sigma^2.
$$
The condition (\ref{con}) is the homogeneity structure of the standard single 
index models for panel data analysis.  $\{G_{1,k}: \ k = 1, \ \cdots, \ H_1\}$ 
and $\{G_{2,k}: \ k = 1, \ \cdots, \ H_2\}$ are unknown partitions.  $H_1$ and 
$H_2$ are unknown integers, $H_1$ is much smaller than $m$, $H_2$ is much 
smaller than $mp$.  $g_{(k)}(\cdot)$, $k=1, \ \cdots, \ H_1$, are unknown 
functions to be estimated, and $\beta_{(k)}$, $k=1, \ \cdots, \ H_2$, are 
unknown parameters to be estimated. 

Let $\beta_{i0}=1$ be the first component of $\vbeta_i$.  In the 
literature, the most commonly used identification condition for the single 
index models is $\|\vbeta_i\|=1$ and $\beta_{i0}>0$, or $\beta_{i0}=1$.  We 
choose the latter in this paper.

The models (\ref{model}) together with (\ref{con}) show that the homogeneity 
pursuit in the single index models for panel data analysis is even more 
important than that in the linear models, this is because we would have to 
estimate $m$ unknown functions and $mp$ unknown parameters in order to explore 
the individual attributes, if the homogeneity pursuit is not conducted.  
However, if the homogeneity pursuit is conducted, we only need to estimate 
$H_1$, much smaller than $m$, unknown functions and $H_2$, much smaller than 
$mp$, unknown parameters when the homogeneity exists.  Even without taking into 
account the benefit resulted from the homogeneity pursuit for the parametric 
part of the models, just for the part of unknown functions alone, to estimate 
much fewer functions would make a big difference in the obtained estimators, in 
terms of the stability of the estimators.

The rest of the paper is organized as follows.  We begin in Section \ref{est} 
with a description of the proposed estimation procedure which is embedded with 
a binary segmentation based homogeneity pursuit.  The asymptotic properties of 
the proposed estimators are presented in Section \ref{asy}.  The performance of 
the proposed estimation procedure and homogeneity pursuit method, when sample 
size is finite, are assessed by simulation studies in Section \ref{simu}.  In 
Section \ref{re}, applying the single index models (\ref{model}) together with 
the homogeneity structure (\ref{con}) to the 49 Industry Portfolios data set, 
which can be freely downloaded from Kenneth French's website 
\begin{center}
{\tt 
http://mba.tuck.dartmouth.edu/pages/faculty/ken.french/data\_library.html}, 
\end{center}
and the UK climate data, which can be freely downloaded from
\begin{center}
{\tt http://www.metoffice.gov.uk/public/weather/climate-historic},
\end{center}
we will show the advantages of the proposed statistical methodology.  We leave 
all technical proofs of the asymptotic properties in the Appendix.

\section{Estimation procedure}
\label{est}

\subsection{Estimation method}
\label{est1}

Our approach to deal with the unknown functions $g_i(\cdot)$, 
$i=1, \ \cdots, \ m$, in (\ref{model}) is based on the B-Spline.  To achieve 
the best result 
for the homogeneity pursuit, we have to decompose all $g_i(\cdot)$s by the same 
B-Spline basis, $\vB(\cdot) = (B_1(\cdot), \ \cdots, \ B_K(\cdot))^{\T}$.

For each $i$, $i = 1, \ \cdots, \ m$, let $\tilde{\vbeta}_i$ be the estimate 
of $\vbeta_i$ obtained, based on the observations for the $i$th individual, 
by a standard estimation procedure for the single index models, e.g. the 
method in \cite{yuruppert02} or in \cite{hardle1989investigating}, and 
$$
a 
= 
\min\limits_{1 \leq i \leq m} \min\limits_{1 \leq t \leq T} 
\vX_{it}^{\T} \tilde{\vbeta}_i,
\quad
b
=
\max\limits_{1 \leq i \leq m} \max\limits_{1 \leq t \leq T}
\vX_{it}^{\T} \tilde{\vbeta}_i.
$$
We use the B-Spline basis of order $s$ in this paper, and the basis, 
$\vB(\cdot)$, is formed by the equally spaced knots, $\tau_k$, 
$k=0, \ \cdots, \ K-s+1$, on the interval $[a, \ b]$, with $\tau_0 = a$ and 
$\tau_{K-s+1} = b$.  Based on the basis $\vB(\cdot)$, $g_i(\cdot)$ can be 
decomposed as
\begin{equation}
g_i(\cdot) \approx \vB(\cdot)^{\T} \vtheta_i,
\label{app}
\end{equation}
where $\vtheta_i = (\theta_{i1}, \ \cdots, \ \theta_{iK})^{\T}$.  So, to get 
the estimator of $g_i(\cdot)$, we only need to get the estimator of $\vtheta_i$.

Our estimation procedure for $\vtheta_i$ and $\vbeta_i$, $i=1, \ \cdots, \ m$, 
consists of three stages: in the first stage, for each $i$, we estimate 
$\vtheta_i$ and $\vbeta_i$ only based
on the observations for the $i$th individual, and treat the obtained
estimators as initial estimators; we identify, in the second stage, the
homogeneity structure in the $\vtheta_i$s and $\vbeta_i$s based on the initial
estimators obtained in the first stage; in the final stage, we estimate the
$\vtheta_i$s and $\vbeta_i$s under the identified homogeneity structure.

We now present the details of the estimation procedure.
\begin{enumerate}
\item[{\bf Stage 1}] ({\it Initial Estimation}).  Let 
$\wb\vbeta_i=(\beta_{i1}, \ \cdots, \ \beta_{ip})\trans$, which 
is $\vbeta_i$ with the first component, which is always $1$, being dropped.  
For each $i$, based on the observations for the $i$th individual, approximating 
$g_i(\cdot)$ by its decomposition (\ref{app}) and applying the least squares 
estimation method, we have the following objective function
\begin{equation}
\sum_{t=1}^T 
\Big(
y_{it} - \vB\trans(\vX_{it}\trans\vbeta_i)\vtheta_i 
\Big)^2.
\label{eqn:min1}
\end{equation}
Minimise (\ref{eqn:min1}) with respect to 
$(\wb\vbeta_i^{\T}, \ \vtheta_i^{\T})$, and denote the resulting minimiser by 
$(\tilde{\vbeta}_i^{\T}, \ \tilde{\vtheta}_i^{\T})$.  We will show how to 
conduct the minimisation in Section \ref{est2}.

\item[{\bf Stage 2}] ({\it Homogeneity Pursuit}).  Let $\tilde{\beta}_{ij}$ be 
the $j$th component of 
$\tilde{\vbeta}_i$, we sort $\tilde{\beta}_{ij}$, $i=1, \ \cdots, \ m$, 
$j=1, \ \cdots, \ p$, in ascending order, and denote them by
$$
b_{(1)}
\leq
\cdots
\leq
b_{(mp)}.
$$
We use $R_{ij}$ to denote the rank of $\tilde{\beta}_{ij}$.  Identifying the 
homogeneity among $\tilde{\beta}_{ij}$, $i=1, \ \cdots, \ m$, 
$j=1, \ \cdots, \ p$, is equivalent to detecting the change points among
$b_{(l)}$, $l=1, \ \cdots, \ mp$.  To this end, we apply the Binary 
Segmentation algorithm as follows. 

For any $1 \leq i <  j \leq mp$, let
$$
\Delta_{ij}(\kappa)
=
\sqrt{\frac{(j-\kappa)(\kappa-i+1)}{j-i+1}}
\left|
\frac{\sum_{l=\kappa+1}^j b_{(l)}}{j-\kappa}
-
\frac{\sum_{l=i}^{\kappa} b_{(l)}}{\kappa-i+1}
\right|.
$$
Given a threshold $\delta$, the Binary Segmentation algorithm to detect the 
change points works as follows
\begin{enumerate}
\item[(1)] Find $\hat{k}_1$ such that
$$
\Delta_{1,mp}(\hat{k}_1)
=
\max\limits_{1 \leq \kappa < mp} \Delta_{1,mp}(\kappa).
$$
If $\Delta_{1,mp}(\hat{k}_1) \leq \delta$, there is no change point among 
$b_{(l)}$, $l=1, \ \cdots, \ mp$, and the process of detection ends.   
Otherwise, add $\hat{k}_1$ to the set of change points and divide the region
$\{\kappa:  \ 1 \leq \kappa \leq mp \}$ into two subregions:
$\{\kappa:  \ 1 \leq \kappa \leq \hat{k}_1 \}$ and
$\{\kappa:  \ \hat{k}_1+1 \leq \kappa \leq mp \}$.

\smallskip
\item[(2)] Detect the change points in the two subregions obtained in (1),
respectively.  Let us deal with the region $\{\kappa:  \ 1 \leq \kappa \leq 
\hat{k}_1 \}$ first.  Find $\hat{k}_2$ such that
$$
\Delta_{1,\hat{k}_1}(\hat{k}_2)
=
\max\limits_{1 \leq \kappa < \hat{k}_1} \Delta_{1,\hat{k}_1}(\kappa).
$$
If $\Delta_{1,\hat{k}_1}(\hat{k}_2) \leq \delta$, there is no change point in 
the region $\{\kappa:  \ 1 \leq \kappa \leq \hat{k}_1 \}$.  Otherwise, 
add $\hat{k}_2$ to the set of change points and divide the region
$\{\kappa:  \ 1 \leq \kappa \leq \hat{k}_1 \}$ into two subregions:
$\{\kappa:  \ 1 \leq \kappa \leq \hat{k}_2 \}$ and
$\{\kappa:  \ \hat{k}_2+1 \leq \kappa \leq \hat{k}_1 \}$.  For the region
$\{\kappa:  \ \hat{k}_1+1 \leq \kappa \leq mp \}$, we find $\hat{k}_3$ such 
that
$$
\Delta_{\hat{k}_1+1,mp}(\hat{k}_3)
=
\max\limits_{\hat{k}_1+1 \leq \kappa < mp} \Delta_{\hat{k}_1+1,mp}(\kappa).
$$
If $\Delta_{\hat{k}_1+1,mp}(\hat{k}_3) \leq \delta$, there is no change point in the region
$\{\kappa:  \ \hat{k}_1+1 \leq \kappa \leq mp \}$.  Otherwise, add $\hat{k}_3$ 
to the set of change points and divide the region $\{\kappa:  \ \hat{k}_1+1 
\leq \kappa \leq mp \}$ into two subregions:
$\{\kappa:  \ \hat{k}_1+1 \leq \kappa \leq \hat{k}_3 \}$ and
$\{\kappa:  \ \hat{k}_3+1 \leq \kappa \leq mp \}$.
\smallskip
\item[(3)] For each subregion obtained in (2), we do exactly the same as that
for the subregion $\{\kappa:  \ 1 \leq \kappa \leq \hat{k}_1 \}$ or
$\{\kappa:  \ \hat{k}_1+1 \leq \kappa \leq mp \}$ in (2), and
keep doing so until there is no subregion containing any change point.
\end{enumerate}

We sort the estimated change point locations in ascending order and denote them 
by
$$
\hat{k}_{(1)} < \hat{k}_{(2)} < \cdots < \hat{k}_{(\hat{H}_{-1})},
$$
where $\hat{H}_{-1}$ is the number of change points detected.  In addition, we 
denote $\hat{k}_{(0)} = 0$, $\hat{H}_2 = \hat{H}_{-1} + 1$, and 
$\hat{k}_{(\hat{H}_2)} = mp$.

We use $\hat{H}_2$ to estimate $H_2$.  Let 
$$
\hat{G}_{2, s} 
= 
\{(i, j): \ \hat{k}_{(s-1)} < R_{ij} \leq \hat{k}_{(s)}\},
\quad
1 \leq s \leq \hat{H}_2,
$$
we use $\left\{\hat{G}_{2, s}: \ 1 \leq s \leq \hat{H}_2 \right\}$
to estimate the partition $\{G_{2,s}: \ 1 \leq s \leq H_2\}$.  We consider all 
the $\beta_{ij}$s with the subscript $(i,j)$ in the same member of the 
estimated partition having the same value.

Let $\tilde{\theta}_{ij}$ be the $j$th component of $\tilde{\vtheta}_i$.  Doing 
exactly the same to $\tilde{\theta}_{ij}$, $i=1, \ \cdots, \ m$, $j=1, \ 
\cdots, \ K$, we get a partition $\{\hat G_{1,1}, \ \cdots, \ \hat G_{1,\hat{H}_1}\}$ of 
$\{(i, j): \ i=1, \ \cdots, \ m; \ j=1, \ \cdots, \ K\}$.  We consider all the  
$\theta_{ij}$s with subscript $(i,j)$ in the same member of the estimated 
partition having the same value.

\item[{\bf Stage 3}] ({\it Final Estimation}).  Let 
$L(\eta_1, \ \cdots, \ \eta_{\hat{H}_2}, \ \xi_1, \ \cdots, \ \xi_{\hat H_1})$
be
\begin{equation}
\label{eqn:min2}
\sum\limits_{i=1}^m \sum_{t=1}^T
\Big(
y_{it} - \vB\trans(\vX_{it}\trans\vbeta_i)\vtheta_i
\Big)^2.
\end{equation}
with $\beta_{ij}$, $i=1, \ \cdots, \ m$, $j=1, \ \cdots, \ p$, being replaced 
by $\eta_k$ if $(i, j) \in \hat{G}_{2, k}$, and $\theta_{ij}$, $i=1, \ \cdots, 
\ m$, $j=1, \ \cdots, \ K$, being replaced by $\xi_s$ if 
$(i, j) \in \hat G_{1,s}$.  Let 
$(\hat{\eta}_1, \ \cdots, \ \hat{\eta}_{\hat{H}_2}, \ \hat{\xi}_1, \ \cdots, \ 
\hat{\xi}_{\hat H_1})$ minimise 
$L(\eta_1, \ \cdots, \ \eta_{\hat{H}_2}, \ \xi_1, \ \cdots, \ \xi_{\hat H_1})$.  
The final estimator $\hat{\beta}_{ij}$ of $\beta_{ij}$ is $\hat{\eta}_k$ if 
$(i, j) \in \hat{G}_{2, k}$, and the final estimator $\hat{\theta}_{ij}$ of 
$\theta_{ij}$ is $\hat{\xi}_s$ if $(i, j) \in \hat G_{1,s}$.  Once we have the 
estimator $\hat{\theta}_{ij}$, the estimator $\hat{g}_i(\cdot)$ of $g_i(\cdot)$ 
is taken to be $\vB(\cdot)^{\T} \hat{\vtheta}_i$.

\end{enumerate}

\bigskip

\noindent
\begin{rmk}  When dealing with the unknown functions $g_i(\cdot)$,
$i=1, \ \cdots, \ m$, in the estimation procedure, instead of treating each 
unknown function as a single undivided unit to conduct homogeneity pursuit, we 
work on the coefficients of its B-Spline decomposition.  This is because there 
may still be some kind of homogeneity between two functions even if they are 
different.  For example, for two different functions, it could be the case that 
some coefficients of the B-Spline decomposition of one function are the same as 
some coefficients of the B-Spline decomposition of another one.  If we treat 
each unknown function as a single undivided unit to conduct homogeneity 
pursuit, we would not identify or use this kind of homogeneity, which would 
make our final estimators not as efficient as they should. 
\end{rmk}

\subsection{Computational algorithm}
\label{est2}

In the estimation procedure described in Section \ref{est1}, the minimiser of 
(\ref{eqn:min1}) does not have a closed form, neither does the minimiser of 
$L(\eta_1, \ \cdots, \ \eta_{\hat{H}_2}, \ \xi_1, \ \cdots, \ \xi_{\hat H_1})$.  To 
conduct the minimisation of either of the two objective functions, we appeal to
the standard NLS algorithm, and use the {\tt nlsLM of minpack.lm} package in R 
to implement it.   One 
can also use other NLS software, for example, the NLS routine {\tt lsqnonlin()} 
from MATLAB and {\tt PROC NLIN} from SAS.  To use the {\tt nlsLM of minpack.lm} 
package in R, we first need to find an initial value.  The initial value for 
minimising (\ref{eqn:min1}) can be obtained as follows:
\begin{enumerate}
\item[(1)] Apply the standard least squares estimation for the linear models to 
$(y_{it}, \ \vX_{it})$, $t = 1, \ \cdots, \ T$, and denote the resulting 
estimator by $\check{\vbeta}_i$, the initial value for $\vbeta_i$ is taken to 
be $\vbeta_i^{(0)} = \check{\beta}_{i0}^{-1} \check{\vbeta}_i$, 
$\check{\beta}_{i0}$ is the first component of $\check{\vbeta}_i$.

\item[(2)] Substitute $\vbeta_i^{(0)}$ for $\vbeta_i$ in (\ref{eqn:min1}), 
then minimise (\ref{eqn:min1}) with respect to $\vtheta_i$, the minimiser 
$\vtheta_i^{(0)}$ is the initial value of $\vtheta_i$.
\end{enumerate}
Once we have $\vbeta_i^{(0)}$ and $\vtheta_i^{(0)}$, the minimiser of 
(\ref{eqn:min1}) can be obtained by the {\tt nlsLM of minpack.lm} package in R 
straightforwardly.

For any set $A$, let $|A|$ be the number of elements in $A$.  The initial 
value for minimising $L(\eta_1, \ \cdots, \ \eta_{\hat{H}_2}, \ \xi_1, \ 
\cdots, \ \xi_{\hat H_1})$ can 
be obtained through the initial estimates of $\vbeta_i$ and $\vtheta_i$, 
obtained in Stage 1 of the estimation procedure in Section \ref{est1}, as 
follows:
$$
\eta_s^{(0)} 
= 
\left(|\hat{G}_{2,s}| \right)^{-1}
\sum\limits_{(i,j) \in \hat{G}_{2,s}} \tilde{\beta}_{ij},
\quad
s = 1, \ \cdots, \ \hat{H}_2
$$
and
$$
\xi_s^{(0)}
=
\left(|\hat G_{1,s}| \right)^{-1}
\sum\limits_{(i,j) \in \hat G_{1,s}} \tilde{\theta}_{ij},
\quad
s = 1, \ \cdots, \ \hat H_1.
$$
Once we have the initial value $(\eta_1^{(0)}, \ \cdots, \ 
\eta_{\hat{H}_2}^{(0)}, \ \xi_1^{(0)}, \ \cdots, \ \xi_{\hat H_1}^{(0)})$, we 
can have the minimiser of $L(\eta_1, \ \cdots, \ \eta_{\hat{H}_2}, \ 
\xi_1, \ \cdots, \ \xi_{\hat H_1})$ by using the {\tt nlsLM of minpack.lm} 
package in R straightforwardly.

\subsection{Selection of tuning parameters}
\label{tune}

The threshold $\delta$ in the Stage 2 of the proposed estimation procedure, 
described in Section \ref{est1}, plays a key role for the success of the 
homogeneity pursuit.  As far as the implementation of the homogeneity pursuit 
is concerned, the selection of $\delta$ is equivalent to the selection of 
$\hat H_1$ and $\hat H_2$, and to select an integer is easier, therefore, 
in this section, instead of selecting $\delta$, we develop a cross-validation 
procedure to select the two tuning parameters, $\hat H_1$ and $\hat H_2$.

For the single index model (\ref{model}) where $\bX_{it}$'s are independent 
across $t=1, \ \cdots, \ T,$  we implement a $L$-fold cross validation 
approach.  In particular, for a given pair $\{H_1, \ H_2\}$, we remove 
$1/L$th of the observed time points for 
$\{(y_{it}, \ \vX_{it}), \ i=1, \ \cdots, \ m, \ t=1, \ \cdots, \ T\}$ 
as a validation set, estimate the single index model (\ref{model}) with 
identified homogeneity structure on the remaining data, compute the squared 
error between $y_{it}$ and fitted values 
$\hat g_i(\vX_{it}^{\T}\hat\vbeta_i)
=\vB\trans(\vX_{it}\trans\hat\vbeta_i)\hat\vtheta_i$, on the validation set, 
and repeat this procedure $L$ times to calculate the cross-validated mean 
squared error and its corresponding standard error.  We search over a grid of 
$\{H_1, \ H_2\}$ values and apply the one-standard-error rule to choose the 
smallest model for which the estimated cross-validated error is within one 
standard error of the lowest point on the error surface.  The rationale here is 
that if a set of models appear to be more or less equally good, then we might 
tend to choose the simplest model.  Across the candidate pairs, 
$\{\hat H_1, \ \hat H_2\},$ whose corresponding errors are within this 
deviation, one can choose the smallest $\hat H_1$ after selecting the smallest 
$\hat H_2$ or switch the selection order or select the smallest value of 
$\hat H_1+\hat H_2$, we take the first approach since it produces better model 
selection consistency in our numerical experiments.  A similar 
one-standard-deviation-rule technique has been adopted to choose the 
regularisation parameter with a smaller model size for the lasso problems 
\cite[]{james13}.

When $\bX_{it}$'s are time dependent panel data, we implement a rolling 
procedure to perform cross-validation for time series.  More specifically, 
for each $r=L, \ L-1, \ \cdots, \ 1,$ we rollingly treat 
$\{(y_{it}, \ \vX_{it}), \ i=1, \ \cdots, \ m, \ t=1, \ \cdots, \ T-r\}$ as 
training observations and $\{(y_{i,T-r+1}, \ \vX_{i,T-r+1}), \ i=1, \ \cdots, \ 
m\}$ as validation set, calculate the squared error between each $y_{it}$ and 
its fitted value.  Finally, we apply the one-standard-deviation-rule on the 
lowest cross-validated mean squared error and choose $\hat H_2$ and $\hat H_1.$

In the cross-validation procedure when we need make predictions for validation 
set, the domain in $\vB(\vX_{it}^{\T}\hat \vbeta_i)$ for traning data set might 
not cover that for validation set.  We adopt the idea in \cite{wangyang09} by 
mapping $\vX^{\T}_{it}\vbeta_i$ to 
$F_i(\vX^{\T}_{it}\vbeta_i) \sim \text{Unif}[0, \ 1]$, 
where $F_i$ is the distribution function of $\vX^{\T}_{it}\vbeta_i$.  We then 
implement the estimation procedure described in Section~\ref{est1} by 
decomposing 
$g_i(\vX^{\T}_{it}\vbeta_i) 
\approx 
\vB^{\T}\big(F_i(\vX^{\T}_{it}\vbeta_i)\big) \vtheta_i$.  
The proposed approach is thus able to make predictions and, as demonstrated by 
some numerical studies, provides very similar sample performance in terms of 
estimation accuracy.

\subsection{Post-processing step}
\label{post}
We equip the Binary-Segmentation-algorithm-based homogeneity pursuit with an 
additional step aimed to enhance the accuracy of detected change-points 
locations through a fine-scale search. To be specific, at each change-point, 
we re-calculate  $\Delta_{ij}(\kappa)$ over the interval between two adjacent 
change-points and identify the new change-point location to replace the old 
one. We perform this post-processing procedure by iteratively cycling 
through all neighbouring change-points and fine-tuning the change-points 
locations. This procedure is terminated when the set of change-points does not 
change. Our numerical experiments show that this extra post-processing step 
apparently improve the accuracy of each estimated change-point location and 
hence the identified homogeneity structure for model~\eqref{model}.

\section{Asymptotic properties}
\label{asy}

In this section, we are going to investigate the asymptotic behaviour of the 
estimators obtained by the proposed estimation procedure, which we call 
correct-fitting, and compare with the estimators obtained without homogeneity 
pursuit, which is the initial estimators obtained in the Stage 1 in the 
proposed estimation procedure, we call it over-fitting, and the estimators 
obtained under the assumption that all individuals share the same index 
(namely, $\vbeta_1 = \cdots = \vbeta_m$), which we call under-fitting.  The 
asymptotic theory presented in this section is in the sense that 
$T\longrightarrow \infty$, and $m$, $p$ are all possibly diverging to infinity 
but $H_1$, $H_2$ are fixed.  This agrees with many applications in which $H_1$ 
and $H_2$ are expected to be small and thus significant reduction of unknown 
parameters can be achieved by clustering the parameters.  To make the 
presentation neat, we state the asymptotic theorems in this section and leave 
all technical proofs in the Appendix.

Let $\epsilon_{t}=(\epsilon_{1t}, \ \cdots, \ \epsilon_{mt})\trans$,
$$
\wb\vX_{it}=(X_{it,1}, \ \cdots, \ X_{it,p})\trans,
\quad
y_t=(y_{1t}, \ \cdots, \ y_{mt})\trans, 
\quad
\vX_{t}=(\vX_{1t}\trans, \ \cdots, \ \vX_{mt}\trans)\trans.
$$
In this paper, we assume $(y_{t}, \ \vX_{t}, \ \epsilon_t)$ are stationary 
with $\alpha(l)\le \rho^l$ for some $\rho<1$, and $\epsilon_t$ is 
independent of $\vX_t$.  Note that unlike \cite{vogt15}, we do not need to 
assume independence or stationarity of variables cross $i$.

We start with the asymptotic properties of the estimators obtained without 
homogeneity pursuit. The convergence rate of the estimator $\tilde{\beta}_{ij}$ is of order $T^{-1/2}$, and the convergence rate of the estimator 
$\tilde{g}_i(u)$ is of order $T^{-2/5}$, which is as expected as we assumed the 
functions are twice differentiable.

\bigskip

\begin{thm}\label{thm:new1}
({\it Over-fitting case}).  For any $i$, $i=1, \ \cdots, \ m$, and 
$1\le j\le p$, under the conditions (C1)-(C4) and (C5') in the Appendix, we 
have
$$
T^{1/2} (\wt\ve_{ij}\trans\wt\vTheta_2\wt\ve_{ij})^{-1/2}
\left(\tilde{\beta}_{ij} - \beta_{ij} \right) 
\convD N(0,\ 1)
$$
and
$$
T^{2/5}(\wt\vb_i\trans(u)\wt\vTheta_1\wt\vb_i(u))^{-1/2}
\left(\wt{g}_i(u) - g_{i}(u)-r_i(u) \right) 
\convD N(0,\ 1),
$$
where $\wt\ve_{ij}$ and $\wt\vb_i(u)$ are unit vectors,  $\wt\vTheta_1$, 
$\wt\vTheta_2$ are matrices with eigenvalues bounded and bounded away from 
zero, all these quantities are defined in the proof in the 
Appendix \ref{sec:conv4}.  The bias term 
$r_i(u)=g_{i}(u)-\vB\trans(u)\vtheta_{0i}$ satisfies $|r_i(u)|\le CK^{-2}$, 
where $\vtheta_{0i}$ is the vector of spline coefficients used to approximate 
$g_i$ as defined in Appendix \ref{sec:ass}.
parameter of $g_i$ as defined in assumption  (C3).
\end{thm}

\bigskip

Let $m_i$ be the size of $G_{1,h}$ that contains $i$, and $m_{ij}$ be the size 
of $G_{2,h}$ that contains $\beta_{ij}$.  To make the statement about the 
correct-fitting case cleaner, we assume that all $m_i$ are of the same order 
and all $m_{ij}$ are of the same order ($\max_{i,j} m_{ij}/\min_{i,j} m_{ij}$ 
and $\max_i m_i/\min_i m_i$ are bounded) in the following theorem, which shows 
in particular that the convergence rate of the estimator $\hat{\beta}_{ij}$ is 
of order $(mpT)^{-1/2}$, and the convergence rate of the estimator 
$\hat{g}_i(u)$ is of order $(mT)^{-2/5}$.

\bigskip

\begin{thm}\label{thm:new2} ({\it Correct-fitting case}).  For any $i$, 
$i=1, \ \cdots, \ m$, and $1\le j\le p$, under the conditions (C1)-(C6) in the 
Appendix, we have
$$
(mpT)^{1/2} (\ve_{ij}\trans\vTheta_2\ve_{ij})^{-1/2} 
\left(\hat{\beta}_{ij} - \beta_{ij} \right) 
\convD N(0, \ 1)
$$
and
$$
(mT)^{2/5}(\vb_i\trans(u)\vTheta_1\vb_i(u))^{-1/2}
\left(\hat{g}_i(u) - g_i(u)-r_i(u) \right) 
\convD N(0, \ 1),
$$
where $\ve_{ij}$ and $\vb_i(u)$ are unit vectors, $\vTheta_1$, $\vTheta_2$ are 
matrices with eigenvalues bounded and bounded away from zero, all these 
quantities are defined in the proof in the Appendix \ref{sec:conv4}. 
\end{thm}

\bigskip

Finally, for the under-fitting case, let $\check{\vbeta}_i$ and 
$\check{g}_i(\cdot)$ be the estimators of $\vbeta_i$ and $g_i(\cdot)$ obtained 
under the assumption that all individuals share the same unknown parameters.

\bigskip

\begin{thm}\label{thm:new3} ({\it Under-fitting case}).  Suppose the 
$\vbeta_i$s are sufficiently separated in the sense that for 
$\bar\vbeta:=\sum_{i=1}^m\vbeta_i/m$, 
$$\frac{1}{mp}\sum_{i=1}^m\|\vbeta_i-\bar\vbeta\|^2\ge c$$
for some $c>0$, then 
$$\frac{1}{mp}\sum_{i=1}^m\|\check \vbeta_i-\vbeta_i\|^2\ge c.$$
Similarly, if 
$$\frac{1}{m}\sum_{i=1}^m\int\left |{g}_i(u) - \bar g(u) \right|^2 du \ge c, 
$$
where $\bar g(u)=m^{-1}\sum_{i=1}^m g_i(u)$, then
$$\frac{1}{m}\sum_{i=1}^m\int\left |\check{g}_i(u) -  g_i(u) \right|^2 du \ge c.
$$
\end{thm}

\section{Simulation studies}
\label{simu}

In this section, we are going to use a simulated example to demonstrate how 
accurate the proposed estimation is.  We will also show much loss it would 
inflict if the homogeneity structure is ignored or mistakenly specified as that 
all individuals share the same index coefficients or the same link function.

\noindent
{\bf Example}.  We generate a sample from model (\ref{model}) with $p=2$ and 
WLOG an even $m,$ 
$$
g_i(u) 
= 
\begin{cases} 
\sin(\pi u/4) & \text{when } i=1, \ 2, \ \cdots, \ m/2, 
\\ 
\cos(\pi u/4) & \text{when } i=m/2+1, \ \cdots, \ m, 
\end{cases}
$$
and 
$$
\vbeta_i 
= 
\begin{cases} 
(1, \ -1.5\sqrt{0.2}, \ -0.5\sqrt{0.2})^{\T} & 
\text{when } i=1, \ 3, \ \cdots, \ m-1,
\\
(1, \ 0.5\sqrt{0.2}, \ 1.5\sqrt{0.2})^{\T} & 
\text{when } i=2, \ 4, \ \cdots, \ m, 
\end{cases}
$$
where $\|\vbeta_i\|^2=1.5$ for $i=1, \ \cdots, \ m$.  Let $\vX_{it}$ and 
$\epsilon_{it}$, $i=1, \ \cdots, \ m$, $t=1, \ \cdots, \ T$ be independently 
generated from $\frac{1}{\sqrt{1.5}}N(\0_3, \ \vI_3)$ truncated by 
$[-1.343, \ 1.343]^3$ (the range of 5th to 95th quantiles for $N(0, \ 2/3)$) 
and $N(0, \ \sigma^2)$, respectively.  Once $\vX_{it}$ and $\epsilon_{it}$ are 
generated, $y_{it}$ can be generated through (\ref{model}).

We conduct the simulated example for various $m$s and $T$s with $\sigma=0.2$, 
and compare our proposed approach to its potential competitors based on the 
following performance metrics: 
\begin{enumerate}
\item[(1)]{\it Estimation accuracy}.  For an estimator $\hat{\vbeta}_i$ of 
$\vbeta_i$, we use the mean squared error (MSE), namely 
$\MSE(\hat{\vbeta}_i) = E\left(\|\hat{\vbeta}_i - \vbeta_i \|^2\right)$, to 
assess the estimation 
error of $\hat{\vbeta}_i$.  Analogously, for an estimator $\hat g_i(\cdot)$ of 
$g_i(\cdot)$, its estimation accuracy can be evaluated based on the mean 
integrated squared error, 
$$
\text{MISE}(\hat g_i)
=
E \left\{\int\left(\hat{g}_i(u) - g_i(u)\right)^2 du \right\}. 
$$ 
To avoid the situation where the performance is dominated by the poor boundary 
behaviour, we let the integral domain to be non-boundary region, which is 
between the 1st and 99th quantiles of 
$\{\vX_{it}^{\T}\vbeta_i, \ t=1, \ \cdots, \ T\}$.

\item[(2)]{\it Homogeneity structure identification consistency}.  To evaluate 
the distance between the detected homogeneity structure and the true one, we 
use the normalized mutual information (NMI) \cite[]{ke15}, which measures the 
similarity between two partitions.  Suppose 
$\mathbb{C}=\{C_1, \ C_2, \ \cdots \}$ and 
$\mathbb{D}=\{D_1, \ D_2, \ \cdots \}$ are two partitions of 
$\{1, \ \cdots, \ n\}$, the NMI is defined as 
$$
\text{NMI}(\mathbb{C}, \ \mathbb{D})
=
\frac{I(\mathbb{C}, \ \mathbb{D})}
{[H(\mathbb{C})+ H(\mathbb{D})]/2},
$$ 
where 
$$
I(\mathbb{C}, \ \mathbb{D}) 
= 
\sum\limits_{k,j} 
\big(|C_k \cap D_j|/n\big) \log\big(n| C_k \cap D_j|/|C_k||D_j|\big)
$$ 
and 
$$
H(\mathbb{C})
=
-\sum\limits_{k} 
\big(|C_k|/n\big) \log\big(|C_k|/n\big).
$$ 
The NMI takes values in $[0, \ 1]$ with larger values indicating higher level 
of similarity between two partitions.  For an estimated partition 
$\hat{\mathbb{G}}_2=\{\hat G_{2,1}, \ \cdots, \ \hat G_{2, \hat H_2}\}$ of 
$\{(i, \ j): \ 1, \ \cdots, \ m, \ j=1, \ \cdots, \ p\}$, obtained in the 
Stage 2 of the proposed estimation procedure in Section \ref{est1}, we 
calculate $\text{NMI}(\hat{\mathbb{G}}_2, \ \mathbb G_2)$ to assess how close 
to the true homogeneity structure in $\beta_{ij}$s the estimated one is.  
Similarly, for an estimated partition $\hat{\mathbb{G}}_1$ of 
$\{i: \ 1, \ \cdots, \ m\}$, we use 
$\text{NMI}(\hat{\mathbb{G}}_1, \ \mathbb G_1)$ to evaluate how close 
the estimated homogeneity structure in $g_{i}(\cdot)$s is to the true one.
\end{enumerate}

For each case, we apply either the single index model (\ref{model}) with the 
standard estimation procedure, the initial estimation of the proposed 
estimation procedure in Section~\ref{est1}, which we call over-fitting (Over), 
the single index model (\ref{model}) with the homogeneity structure 
(\ref{con}) together with the proposed estimation procedure, which we call 
correct-fitting, the single index model (\ref{model}) with all individuals 
share the same index vector (namely, $\vbeta_1 = \cdots = \vbeta_m$), 
which we call Under-I, the single index model (\ref{model}) with all 
individuals share the same link function (namely, 
$g_1(\cdot) = \cdots = g_m(\cdot)$, i.e. $\vtheta_1 = \cdots = \vtheta_m$), 
which we call Under-F, or the single index model (\ref{model}) with all 
individuals share both the same index and link function, which we call 
Under-I-F, to the simulated data set.

We develop three methods under the correct-fitting case.  The first approach, 
named Correct-C, optimises (\ref{eqn:min2}) based on the estimated 
{\it componentwise} homogeneity structure in $\beta_{ij}$s and $\theta_{ij}$s, 
obtained in the Stage 2 of the proposed estimation procedure in Section 
\ref{est1} with the tuning parameters selected through the cross-validation 
approach described in Section~\ref{tune}.  The second approach, Correct-V, is 
the same as the first approach but optimises (\ref{eqn:min2}) based on the 
estimated componentwise homogeneity structure in $\beta_{ij}$s and 
{\it vectorwise} homogeneity structure in $\vtheta_i$s which can be obtained 
through the estimated componentwise homogeneity structure in $\theta_{ij}$s.  
The third approach, which we call Correct-NMI, is the same as the second 
approach but with the tuning parameters selected to be the one maximising 
$\text{NMI}(\hat{\mathbb{G}}_2, \ \mathbb G_2)$ and 
$\text{NMI}(\hat{\mathbb{G}}_1, \ \mathbb G_1)$.  In practice without knowing 
the true homogeneity structure, one cannot implement Correct-NMI.

Under-I, Under-F and Under-I-F are three kinds of under-fitting.  For Under-I 
or Under-F, the homogeneity structure in $\vtheta_i$s or $\beta_{ij}$s is 
estimated in the same way as that in the proposed estimation 
procedure in Section \ref{est1} with the tuning parameters still selected by 
the one-standard-deviation-rule cross-validation approach.

We compare over-fitting, correct-fittings and under-fittings to the oracle case 
where the true homogeneity structure is used.  The computational algorithms 
for the under-fitting and oracle estimators are the same as that for the 
correct-fitting, but use either identified or pre-specified homogeneity 
structure.  We compare the sample performance of all eight approaches in our 
conducted simulation study.

We report the results for estimation errors and NMIs for $\vbeta_i$s and 
$g_i(\cdot)$s averaged over 100 replicates in Tables~\ref{tab_beta} and 
\ref{tab_fun}, respectively.  In terms of estimation error, the overall 
estimation accuracy is improved as $m$ and $T$ increase and three 
correct-fitting approaches perform very well as reflected in their lower values 
of MSEs and MISEs.  Among the three methods, Correct-NMI provides the best 
performance even producing very comparable MSEs and MISEs with the oracle 
estimator and Correct-C is outperformed by Correct-V in most settings.  This is 
somewhat expected, since, unlike Correct-C, which optimises (\ref{eqn:min2}) 
based on the detected homogeneity structure in $\beta_{ij}$s and 
$\theta_{ij}$s, Correct-V separates the final estimation step from the 
cross-validation procedure, which is used to identify the homogeneity structure 
in $\beta_{ij}$s and $\vtheta_{i}$s.   Analogously, Correct-NMI solves a 
separate optimisation after detecting the homogeneity structure based on the 
largest NMIs.  It is also worth noting that the over-fitting and under-fitting 
methods, which either ignores or mistakenly specify the homogeneity structure, 
provide much worse results, highlighting the importance of incorporating the 
appropriate homogeneity structure.  In terms of selecting the structure of 
homogeneity, we observe that three correct-fitting methods produce perfect 
identifications of homogeneity structure in $\beta_{ij}$s and Correct-NMI 
provides the largest NMI values indicating that it can effectively recover the 
true homogeneity structure in $g_i(\cdot)$s.  The performance of Correct-C and 
Correct-V deteriorates when $m$ increases, this is intuitively due to the 
increased $m$ values and the cross-validation procedure, which tends to choose 
a larger number of change points as $m$ increases, resulting in smaller NMI 
values for $\hat{\mathbb{G}}_1$.

\begin{table}[htbp]
	\begin{center}
		\caption{{\bf The Average of $\MSE(\hat{\vbeta}_i)$, 
$i=1, \ \cdots, \ m$, and Average NMIs for $\hat{\mathbb{G}}_2$} 
\centerline{{\it All entries for MSEs are $10^4$ times their actual values}}}
		\label{tab_beta}
		\vspace{0.3cm}
		\begin{tabular}{cc|ccc|ccc}
			\hline 
			\hline
			& $T$ & \multicolumn{3}{c|}{400} & \multicolumn{3}{c}{800}\tabularnewline
			& $m$ & 30 & 60 & 90 & 30 & 60 & 90\tabularnewline
			\hline 
			\multirow{8}{*}{MSE} 
			& Oracle & 0.380 & 0.205 & 0.133 & 0.214 & 0.128 & 0.101\tabularnewline
			& Correct-C & 0.473 & 0.221 & 0.147 & 0.157 & 0.093 & 0.085 \tabularnewline
			& Correct-V & 0.381 & 0.208 & 0.134 & 0.214 & 0.129 & 0.100 \tabularnewline
			& Correct-NMI & 0.381 & 0.206 & 0.133 & 0.214 & 0.128 & 0.101 \tabularnewline
			& Over & 5.438 & 5.285 & 5.246 & 2.636 & 2.689 & 2.663 \tabularnewline
			& Under-I-F & 4005.6 & 4002.6 & 4002.1 & 4003.4 & 4001.5 & 4001.1\tabularnewline
			& Under-I & 4002.4 & 4001.7 & 4001.3 & 4001.8 & 4001.0 & 4000.6\tabularnewline
			& Under-F & 21.233 & 3.393 & 2.259 & 3.111 & 1.454 & 1.123\tabularnewline
			\hline 
			\multirow{8}{*}{NMI} 
			& Oracle & 1.000 & 1.000 & 1.000 & 1.000 & 1.000 & 1.000\tabularnewline
			& Correct-C & 1.000 & 1.000 & 1.000 & 1.000 & 1.000 & 1.000\tabularnewline
			& Correct-V & 1.000 & 1.000 & 1.000 & 1.000 & 1.000 & 1.000\tabularnewline
			& Correct-NMI & 1.000 & 1.000 & 1.000 & 1.000 & 1.000 & 1.000\tabularnewline
			& Over & 0.339 & 0.290 & 0.267 & 0.339 & 0.290 & 0.267\tabularnewline
			& Under-I-F & 0 & 0 & 0 & 0 & 0 & 0\tabularnewline
			& Under-I & 0 & 0 & 0 & 0 & 0 & 0\tabularnewline
			& Under-F & 0.817  & 0.810 & 0.807 & 0.812 & 0.809 & 0.806\tabularnewline
			\hline 
			\hline
		\end{tabular}
	\end{center}
\end{table}

\begin{table}[htbp]
	\begin{center}
		\caption{{\bf The Average of $\MISE(\hat{g}_i)$, $i=1, \ 
\cdots, \ m$, and Average NMIs for $\hat{\mathbb{G}}_1$} 
\centerline{{\it 
where all entries for MISEs are $10^2$ times their actual values
}}}
		\label{tab_fun}
				\vspace{0.3cm}
		\begin{tabular}{cc|ccc|ccc} 
			\hline
			\hline
			& $T$ & \multicolumn{3}{c|}{400} & \multicolumn{3}{c}{800}\tabularnewline
			& $m$ & 30 & 60 & 90 & 30 & 60 & 90\tabularnewline
			\hline 
			\multirow{8}{*}{MISE} 
			& Oracle & 0.260 & 0.249 & 0.243 & 0.251 & 0.241 & 0.236\tabularnewline
			& Correct-C & 0.667 & 0.553 & 0.485 & 0.366 & 0.318 & 0.307 \tabularnewline
			& Correct-V & 0.300 & 0.289 & 0.311 & 0.257 & 0.246 & 0.246\tabularnewline
			& Correct-NMI & 0.266 & 0.254 & 0.249 & 0.252 & 0.241 & 0.238\tabularnewline
			& Over & 0.548 & 0.544 & 0.540 & 0.410 & 0.403 & 0.396\tabularnewline
			& Under-I-F & 87.988 & 87.963 & 87.959 & 88.746 & 88.726 & 88.647\tabularnewline
			& Under-I & 10.099 & 9.853 & 9.952 & 9.770 & 9.719 & 9.711\tabularnewline
			& Under-F & 85.196 & 85.240 & 85.235 & 86.017 & 86.035 & 85.952\tabularnewline
			\hline 
			\multirow{8}{*}{NMI} 
			& Oracle & 1.000 & 1.000 & 1.000 & 1.000 & 1.000 & 1.000 \tabularnewline
			& Correct-C & 0.816 & 0.790 & 0.637 & 0.958 & 0.948 &0.881 \tabularnewline
			& Correct-V & 0.816 & 0.790 & 0.637 & 0.958 & 0.948 & 0.881\tabularnewline
			& Correct-NMI & 0.965 & 0.966 & 0.948 & 0.994 & 0.997 & 0.997\tabularnewline
			& Over & 0.339 & 0.290 & 0.267 & 0.339 & 0.290 & 0.267\tabularnewline
			& Under-I-F & 0 & 0 & 0 & 0 & 0 & 0\tabularnewline
			& Under-I & 0.811 & 0.797 & 0.783 & 0.901 & 0.894 & 0.900\tabularnewline
			& Under-F & 0 & 0 & 0 & 0 & 0 & 0\tabularnewline
			\hline
			\hline
		\end{tabular}
	\end{center}
\end{table}

\section{Real data analysis}
\label{re}

We will illustrate the proposed method with two real data examples in this 
section.

\subsection{Industrial Portfolio's return}
\label{port}

We first study the data set about $m=49$ Industrial Portfolios' daily simple 
return from 1/8/2015 to 31/12/2015.  This data set can be freely downloaded 
from Kenneth French's website

\begin{center}
{\tt http://mba.tuck.dartmouth.edu/pages/faculty/ken.french/data\_library.html}
\end{center}

This data set has been analysed in quite a few literature.  For example, 
\cite{guo2016dynamic} used this data set to demonstrate the performance of a 
newly developed dynamic portfolio allocation.  In this paper, we are going to 
explore the homogeneity structure in this data set by our proposed method.

Let $y_{it}$ be the daily simple return of the $i$th portfolio at the $t$th 
day, $i=1, \ \cdots, \ m, \ t=1, \ \cdots, \ T$, and 
$\vX_{it}=(X_{t1}, \ X_{t2}, \ X_{t3})^{\T}$ be the observation of the 
Fama-French three factors, where $X_{t1}, \ X_{t2}, \ X_{t3}$ respectively 
represent the market (Rm-Rf), size (SMB) and value (HML) factors at the $t$th 
day.

We apply the single index model (\ref{model}) with the unknown homogeneity 
structure (\ref{con}) to fit the data set.  From interpretation point of view, 
the homogeneity structure in the unknown link functions, $g_i(\cdot)$s, where 
each $g_i(\cdot)$ is treated as a single undivided unit would make much more 
sense than the homogeneity structure in the coefficients of the B-Spline 
decompositions of $g_i(\cdot)$s.  Therefore, we use the Correct-V, described 
in Section~\ref{simu}, to identify the homogeneity structure in $\beta_{ij}$s 
or $g_i(\cdot)$s, and estimate the unknown parameters and unknown functions.

In the implementation of the Correct-V, we implement the method in 
Section \ref{est1} with the tuning parameters selected by the cross-validation 
for time series as described in Section~\ref{tune}.  Specifically, we define 
the cross-validated mean squared error
\begin{equation}
\label{cv.err}
\CV
= 
\frac{1}{m L}\sum\limits_{i=1}^{m} \sum\limits_{t=T-L+1}^{T}
\left(y_{it}-\hat{y}_{it}\right)^2
\end{equation}
where $L=30$.  Note that we here do not apply the one-standard-rule when 
performing the cross 
validation to select the tuning parameters for identifying the homogeneity 
structure, since we have already selected a small enough model with 11 and 2 
detected groups in index coefficients and link functions, respectively.  
Table~\ref{tab4} provides the identified clustering results for $\beta_{i2}$, 
$\beta_{i3}$, $g_{i}(\cdot)$, $i=1, \ \cdots, \ 49$, and 
Figure~\ref{funcplot_port} plots the estimated link functions.  We observe a 
few apparent patterns.  Firstly, the estimated link functions are very linear 
indicating the linear relationship between portfolio returns and Fama-French 
three factors, which has been verified by broad empirical studies.  Secondly, 
many portfolios belonging to similar industrials were grouped into the same 
cluster for the estimated index coefficients, e.g. Hardw, Softw and Agric, 
Food, Soda were clustered into Groups~8 and 4 in terms of the estimated 
coefficients for factors SMB and HML respectively.

\begin{table}[htbp]
\footnotesize
	\caption{\label{tab4}{\bf  
Grouping Results for The Index Coefficients for SMB, HML and 
\centerline{Link Functions of 49 Industrial Portfolios}}}
\begin{center}
\begin{tabular}{c|cccccccccc}
	\hline 
	\hline
	& Agric & Food & Soda & Beer & Smoke & Toys & Fun & Books & Hshld & Clths\tabularnewline
	SMB & 6 & 5 & 6 & 4 & 2 & 8 & 9 & 7 & 6 & 6\tabularnewline
	HML & 4 & 4 & 4 & 3 & 4 & 5 & 3 & 4 & 5 & 6\tabularnewline
	Function & ii & i & ii & ii & i & i & i & i & i & ii\tabularnewline
	\hline 
	& Hth & MedEq & Drugs & Chems & Rubbr & Txtls & BldMt & Cnstr & Steel & FabPr\tabularnewline
	SMB & 9 & 8 & 9 & 7 & 7 & 8 & 8 & 8 & 9 & 10\tabularnewline
	HML & 2 & 2 & 1 & 6 & 4 & 4 & 6 & 6 & 9 & 7\tabularnewline
	Function & i & i & i & i & i & i & i & i & i & ii\tabularnewline
	\hline 
	& Mach & ElcEq & Autos & Aero & Ships & Guns & Gold & Mines & Coal & Oil\tabularnewline
	SMB & 8 & 9 & 7 & 6 & 8 & 6 & 10 & 9 & 11 & 9 \tabularnewline
	HML & 8 & 6 & 6 & 6 & 7 & 4 & 10 & 7 & 11 & 9 \tabularnewline
	Function & i & i & i & i & i & i & ii & i & ii & i\tabularnewline
	\hline 
	& Util & Telcm & PerSv & BusSv & Hardw & Softw & Chips & LabEq & Paper & Boxes\tabularnewline
	SMB & 4 & 8 & 8 & 7 & 8 & 8 & 7 & 7 & 6 & 6 \tabularnewline
	HML & 5 & 4 & 5 & 4 & 5 & 3 & 5 & 5 & 5 & 6\tabularnewline
	Function & i & i & i & i & i & i & i & i & i & i \tabularnewline
	\hline 
	& Trans & Whlsl & Rtail & Meals & Banks & Insur & RIEst & Fin & Other & \tabularnewline
	SMB & 7 & 7 & 7 & 7 & 7 & 6 & 7 & 6 & 7 & \tabularnewline
	HML & 6 & 6 & 4 & 3 & 6 & 5 & 5 & 5 & 6 & \tabularnewline
	Function & i & i & i & i & i & i & i & i & ii & \tabularnewline
	\hline 
	\hline
\end{tabular}
\end{center}
\end{table}

\begin{figure*}
	\centering
	\includegraphics[width=16cm,height=16cm]{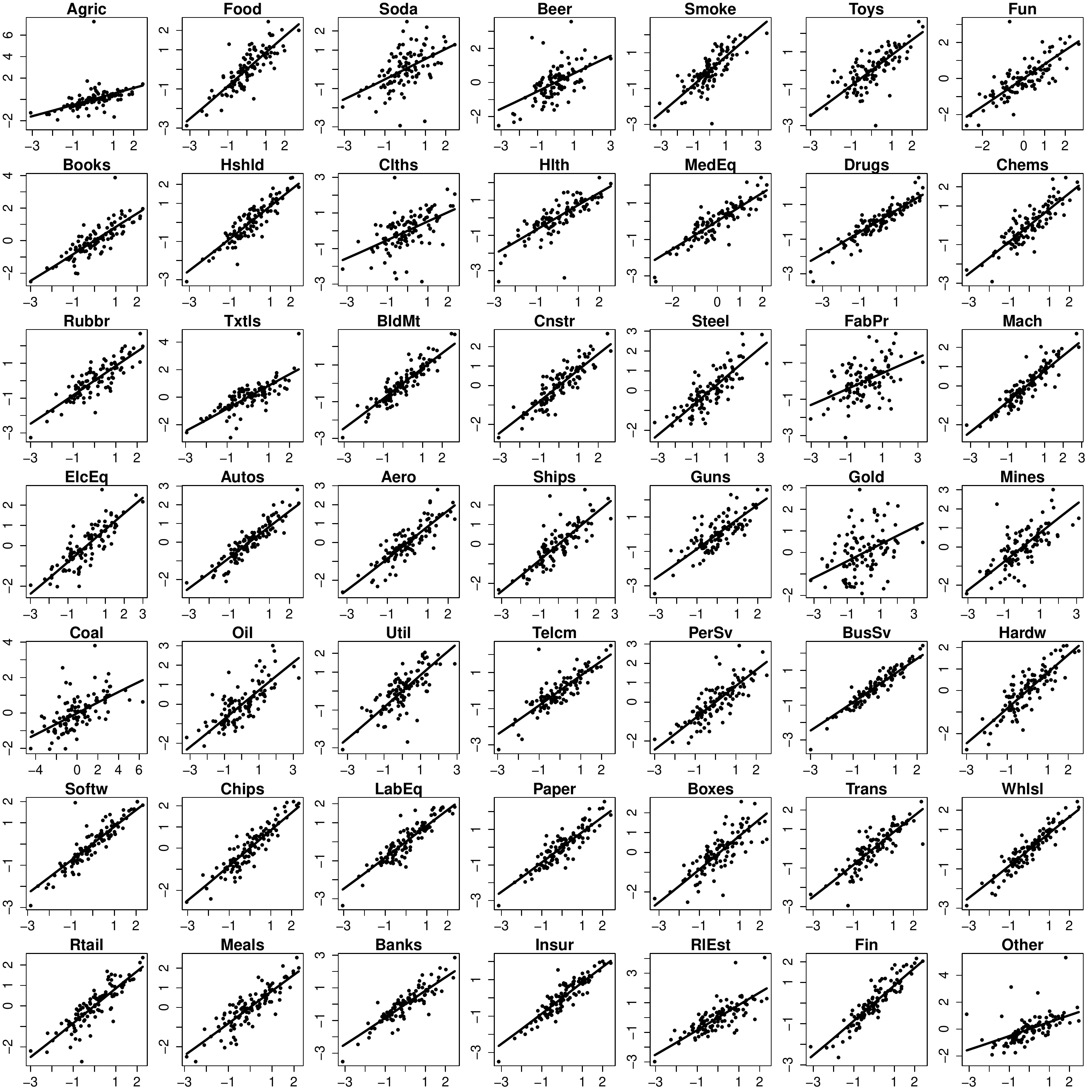}
	\caption{\label{funcplot_port}{Plots of estimated link functions with 
respect to $\{\vX_{it}^{\T}\hat \vbeta_i, t=1,\dots, T\}$.}}
\end{figure*}

\subsection{UK climate data}
\label{temp}

Our second data set, which is available from the UK Met Office website 
\begin{center}
{\tt http://www.metoffice.gov.uk/public/weather/climate-historic}, 
\end{center}
contains monthly data of the mean daily maximum temperature (TMAX), mean daily 
minimum temperature (TMIN), days of air frost (AF), total rainfall (RAIN) and 
total sunshine duration (SUN) collected from 37 stations across the UK.  We 
first remove the missing values and thus select data during the period of 
January 1993 to December 2009 from 16 locations.  We then eliminate the 
seasonality and trend effects and standardise the data.  Let $y_{it}$ be the 
monthly mean temperature, which can be calculated as (TMAX+TMIN)/2, and 
$\vX_{it}=(X_{it1}, \ X_{it2}, \ X_{it3})^{\T}$ be the observations for AF, 
RAIN and SUN, from the $i$th station at the $t$th month, $i=1, \ \cdots, \ 16, 
\ t=1, \ \cdots, \ T$,

Like the analysis of the Industrial Portfolio's return data set, we apply the 
single index model (\ref{model}) with unknown homogeneity structure (\ref{con}) 
together with the proposed estimation procedure, Correct-V, to the data set.  
Table~\ref{tab6} provides the clustering results for the index coefficients and 
link functions, where 4 and 2 groups were selected respectively.  
Figure~\ref{funcplot_temp} plots the estimated link functions at 16 stations.

It is very interesting to see, from Table~\ref{tab6}, that Oxford, Hurn, 
Eastbourne and Bradford share exactly the same model, which implies the impact
of rainfall or total sunshine duration on monthly mean temperature has exactly
the same pattern in these four areas.  The same finding also appears in the 
three areas of Waddington, Sheffield and Heathrow, the two areas of 
Ross-On-Wye and Eskdalemuir, and the two areas of Paisley and Leuchars.  If we 
only focus on the impact of rainfall on monthly mean temperature, the seven 
areas of Waddington, Sheffield, Shawbury, Paisley, Leuchars, Lerwick and 
Heathrow would have exactly the same pattern.  Similar finding also appears 
for the impact of total sunshine duration on monthly mean temperature.

\begin{table}[htbp]
\footnotesize		
	\begin{center}			
	\caption{{\bf Grouping Results for The Index Coefficients for RAIN, 
SUN and 
\centerline{Link Functions at 16 Locations}}}
	\label{tab6}
			\vspace{0.3cm}
		\begin{tabular}{c|cccccccc}
			\hline 
			\hline
			& Waddington & Sheffield & Shawbury & Ross-On-Wye & Paisley & Oxford & Leuchars & Lerwick\tabularnewline
			RAIN & 4 & 4 & 4 & 3 & 4 & 3 & 4 & 4 \tabularnewline
			SUN & 1 & 1 & 2 & 2 & 3 & 2 & 3 & 4 \tabularnewline
			Function & i & i & i & ii & i & i & i & i \tabularnewline
			\hline 
			& Hurn & Heathrow & Eskdalemuir & Eastbourne & Cambridge & Camborne & Bradford & Armagh \tabularnewline
			RAIN & 3 & 4 & 3 & 3 & 3 & 3 & 3 & 4 \tabularnewline
			SUN & 2 & 1 & 2 & 2 & 1 & 3 & 2 & 2 \tabularnewline
			Function & i & i & ii & i & i & i & i & ii \tabularnewline
			\hline 
			\hline
		\end{tabular}
	\end{center}
\end{table}

\begin{figure*}
	\label{fig3}
	\centering
	\includegraphics[width=16cm,height=16cm]{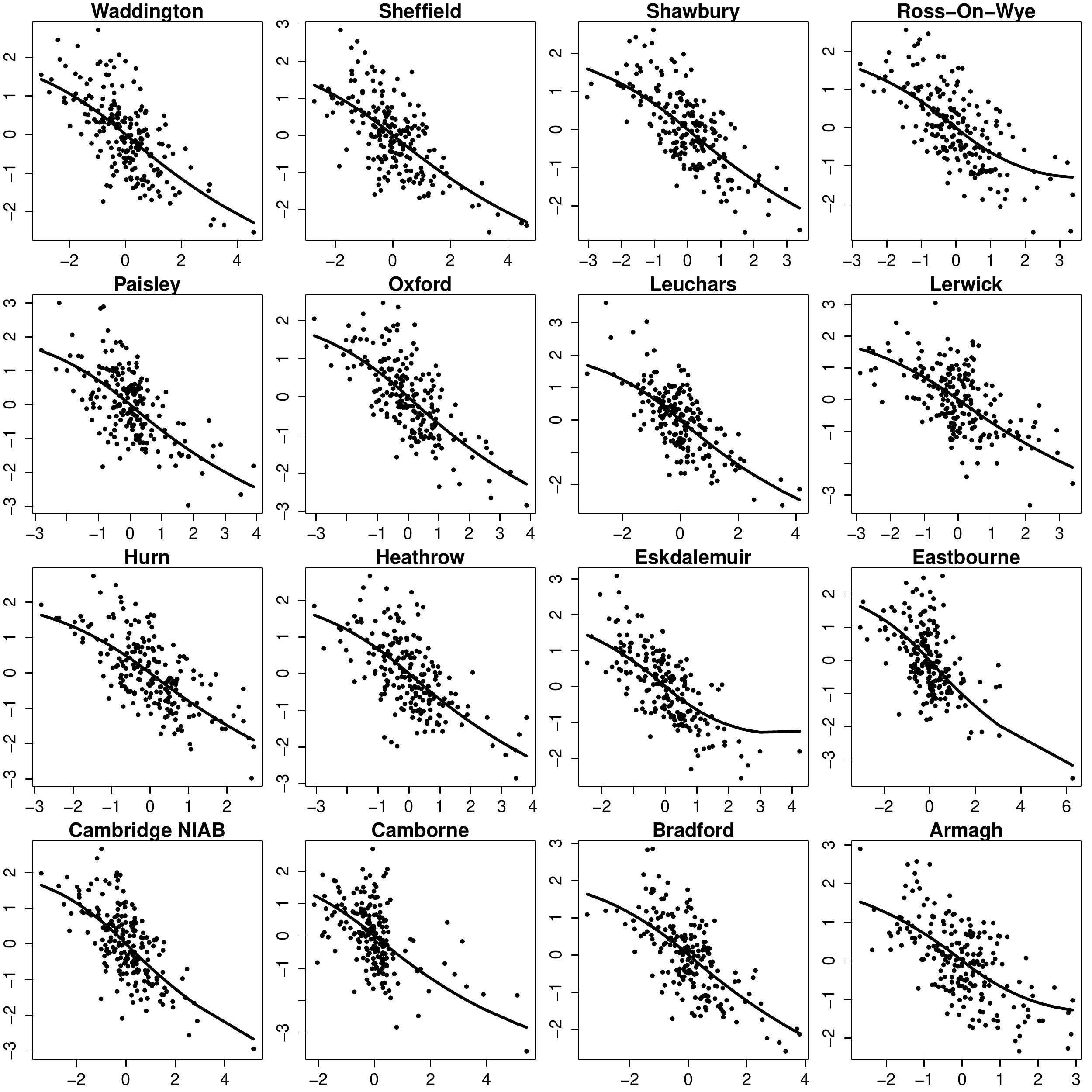}
	\caption{\label{funcplot_temp}{Plots of estimated link functions with 
respect to $\{\vX_i^{\T}\hat \vbeta_{it}, t=1, \dots, T\}$.}}
\end{figure*}

\appendix
\section*{Appendix A. Proofs of Main Results}

\renewcommand{\thesubsection}{A.\arabic{subsection}}

\setcounter{equation}{0}

\renewcommand{\theequation}{A.\arabic{equation}}

\subsection{Assumptions and notations}\label{sec:ass}
Below we use subscript 0 to indicate the true value. We impose  the following assumptions.
\begin{itemize}
\item[(C1)] $(y_t,\vX_t,\epsilon_t), t=1,\ldots,T$ is stationary and $\alpha$-mixing with mixing coefficient $\alpha(l)\le\rho^l$ for some $\rho\in(0,1)$. $\epsilon_{it}$ has mean zero, with variance uniformly bounded, and is independent of $\{\vX_{1t},\ldots,\vX_{mt}\}$. The variables $ X_{it,j}$ are uniformly bounded. The density of $\vX_{it}\trans\vbeta_{0i}$, denoted by $f_i(x)$, is supported on an interval of length, say, $L$ and $Lf_{i}(x)$ is bounded and bounded away from zero on its support, uniformly over $i$.
\item[(C2)] Let $\sigma_{ii',l}=E[\epsilon_{it}\epsilon_{i't'}]$ with $|t-t'|=l$. We assume $\sum_{l=1}^T|\sigma_{ii',l}|\le \tau_{ii'}$ for some $\tau_{ii'}>0$ and $\max_i\sum_{i'}\tau_{ii'}\le M$ for some constant $M$. 
\item[(C3)] The link functions $g_{0i}$ are twice continuously differentiable. We also assume $E[\wb\vX_{it}|\vX_{it}\trans\vbeta_{i}=x]$ is twice continuously differentiable for $\vbeta_i$ in a neighborhood of $\vbeta_{0i}$.
\item[(C4)]   Both $E[\vX_{it}\vX_{it}\trans]$ and $E[(g_{0i}'(\vX_{it}\trans\vbeta_{i}))^2(\wb\vX_{it}-E[\wb\vX_{it}|\vX_{it}\trans\vbeta_{i}])^{\otimes 2}]$ have eigenvalues bounded and bounded away from zero, uniformly over $i$ and $\vbeta_i$ in a neighborhood of $\vbeta_{0i}$, where for any matrix $\vA$, $\vA^{\otimes 2}=\vA\vA\trans$.
\item[(C5)] $H_1$ and $H_2$ are fixed and $\max_{i,j} m_{ij}/\min_{i,j} m_{ij}$ and $\max_i m_i/\min_i m_i$ are bounded, and we set $K\asymp(mT)^{1/5}$. Assume $\frac{m^{6/5}p\log T}{T^{3/5}}\rightarrow 0$ and $\frac{(K^3+p^2)p(\log T)^3}{T}\rightarrow 0$.
\item[(C6)] Assume $mK\sqrt{\log (Tm)/T}<<\delta_1<<\gamma_1$, where $\gamma_1$ is the minimum jump size for the sequence $\theta_{0(1)}\le\cdots\le\theta_{0(mK)}$ at the change points, and $\delta_1$ is the threshold used in the change point detection algorithm (we stop partitioning if the test statistic is below $\delta_1$). Similarly, assume $mp\sqrt{\log (Tm)/T}<<\delta_2<<\gamma_2$, where $\gamma_2$ and $\delta_2$ are similarly defined for the sequence $\beta_{0(1)}\le\cdots\le\beta_{0(mp)}$.

\end{itemize}
\begin{rmk} (C1) contains some mild regularity assumptions. Assuming $X_{it.j}$ to be bounded is common in estimation with B-splines since the basis functions are constructed on a compact interval. If $p$ is fixed, we can simply assume the density of $\vX_{it}\trans\vbeta_{0i}$  is bounded and bounded away from zero. Our assumption however deals with the case $p$ is diverging and thus the length of the support of the density is also diverging. (C2) roughly means the dependence across $i$ is not too strong. If $m$ is fixed, (C2) follows from the geometric mixing assumption. Assumptions similar to (C2) were also used in \cite{bai03} to impose weak dependence among errors. Note \cite{vogt15} made the stronger assumption that the data are independent across $i$ which also easily implies (C2). (C3) contains smoothness condition for some functions and (C4) contains some identifiability conditions usually assumed in single-index models and involves the projection one typically use to profile out the nonparametric part. 
Uniformity over $i$ in various assumptions above is void if $m$ is fixed. (C5) specifies the required divergence rate for $T,m,p,K$. Finally, (C6) is used in showing that  stage 2 of our estimation procedure can identify the true partition with probability approaching one.
\end{rmk}
When considering the estimator in stage 1 of our estimation procedure, we can replace (C5) with the following.
\begin{itemize}
\item[(C5')] We set $K\asymp T^{1/5}$, and assume $(K+p)p\log T/T^{3/5}\rightarrow 0$, $(K^3+p^2)p(\log T)^3/T\rightarrow 0$.
\end{itemize}

Due to assumption (C3), there exists $\vtheta_0=(\vtheta_{01}\trans,\ldots,\vtheta_{0m}\trans)\trans$, $\vtheta_{0i}=(\theta_{0i1},\ldots,\theta_{0iK})\trans$ such that $\sup_x|g_{0i}(x)-\vtheta_{0i}\trans\vB(x)|\le CK^{-2}$. Here and below we use $C$ to denote a generic positive constant whose value can change even on the same line. We use $\|.\|_{op}$ to denote the operator norm of a matrix (the operator norm is the same as the largest singular value) and use $\|.\|$ to denote the Frobenius norm of a matrix. We use $\|.\|_{L^2}$ to denote the $L^2$ norm of functions and $\|.\|_\infty$ is the sup-norm for vectors (maximum absolute value of the components).

Assume the true partition of components of $\vtheta_0$ and $\wb\vbeta_0$ is given by $\cup_{h=1}^{H_1}G_{1,h}=\{1,\ldots,mK\}$ and $\cup_{h=1}^{H_2}G_{2,h}=\{1,\ldots,mp\}$, respectively. The unique values of the components of $\vtheta_0$ and $\wb\vbeta$ are denoted by $\vxi_0=(\xi_{01},\ldots,\xi_{0H_1})\trans\in R^{H_1}$ and $\veta_0=(\eta_{01},\ldots,\eta_{0H_2})\trans\in R^{H_2}$, respectively. Let $\vJ_i^{G_1}$ be the $K\times H_1$ binary matrix whose $(k,h)$ entry is 1 if $\theta_{0ik}=\xi_h$ and 0 otherwise. We  have $\vtheta_{0i}=\vJ_i^{G_1}\vxi_0$. Similarly, we define $\vJ_i^{G_2}$ such that $\vbeta_{0i}=\vJ_i^{G_2}\veta_0$. The sizes of $G_{1,h}$ and $G_{2,h}$ are denoted by $|G_{1,h}|$ and $|G_{2,h}|$, respectively. Finally, let $\vD^{G_1}$ and $\vD^{G_2}$ be the diagonal matrix with entries $\sqrt{|G_{1,h}|}$ and $\sqrt{|G_{2,h}|}$, respectively.
\subsection{Proof summary}
We first define the oracle estimator as the minimizer $(\wh\vtheta,\wh\vbeta)$ of
$$\min_{\vtheta,\wb\vbeta}\sum_{i=1}^m\sum_{t=1}^T(y_{it}-\vB\trans(\vX_{it}\trans
\vbeta_i)\vtheta_i)^2,$$
where $\vbeta_i=(1,\wb\vbeta_i\trans)\trans=(1,\beta_{i1},\ldots,\beta_{ip})\trans$ and $\vtheta_i=(\theta_{i1},\ldots,\theta_{iK})\trans$
with the constraint that components of $\wb\vbeta=(\wb\vbeta_1\trans,\ldots,\wb\vbeta_m\trans)\trans$ in the same partition take the same value and components of $\vtheta=(\vtheta_1\trans,\ldots,\vtheta_m\trans)\trans$ in the same partition take the same value. Here we assume the partition is the true partition, thus the name ``oracle". To make our arguments applicable to over-fitting case, we note that all arguments carry over when the partition used in the oracle estimator is finer than the true partition and thus Theorem \ref{thm:new1} is actually a special case.

In \ref{sec:conv1}-\ref{sec:conv4}, we show that the oracle estimator satisfies the asymptotic normality properties stated in Theorem \ref{thm:new2} (we also obtained convergence rate and asymptotic normality for the entire vector  $\vbeta$ and $\vtheta$, see for example (\ref{eqn:anbeta}) and (\ref{eqn:antheta})). Also, Theorem \ref{thm:new1} follows directly as a special case that each component of $\vtheta$ and $\vbeta$ forms its own group in the partition. Then we show that the change points can be consistently estimated,  and thus the estimator we obtain in stage 3 will be exactly the same as the oracle estimator using the true partition, with probability approaching one, and Theorem \ref{thm:new2} is proved. 

\subsection{Proof of asymptotic property for the oracle estimator}\label{sec:conv1}

In this part we consider the asymptotic property of the oracle estimator, denoted by $(\wh\vtheta,\wh\vbeta)$ in this section, which assumed knowledge of the true partitions. For clarity of presentation, the proof is split into several steps and the proofs of some lemmas were relegated to Appendix B.

\noindent\textit{STEP 1. Prove the convergence rate $\|\wh\vtheta-\vtheta_0\|+\|\wh\vbeta-\vbeta_0\|=O_p(\sqrt{(H_1+H_2)/T}+\sqrt{m}K^{-2})$.} 

In this section, when we use $\vtheta$, we always assume $\vtheta_i=\vJ_i^{G_1}\vxi$ for some $\vxi\in R^{H_1}$ (that is, components of $\vtheta$ are partitioned in the same way as is the true $\vtheta_{0}$). It is easy to see that $\|\vtheta-\vtheta_{0}\|=\|\vD^{G_1}(\vxi-\vxi_0)\|$. Similarly, we always assume $\vbeta_i=\vJ_i^{G_2}\veta$ for some $\veta\in R^{H_2}$ and $\|\vbeta-\vbeta_{0}\|=\|\vD^{G_2}(\veta-\veta_0)\|$.

Define $r_T=\sqrt{(H_1+H_2)/T}+\sqrt{m}K^{-2}$. We only need to show that
\begin{eqnarray*}
&&\inf_{\|\vbeta-\vbeta_0\|^2+\|\vtheta-\vtheta_0\|^2= Lr_T^2} \sum_{i=1}^m\sum_{t=1}^T (y_{it}-\vtheta_i\trans\vB(\vX_{it}\trans\vbeta_i))^2- \sum_{i=1}^m\sum_{t=1}^T (y_{it}-\vtheta_{0i}\trans\vB(\vX_{it}\trans\vbeta_{0i}))^2>0  
\end{eqnarray*}
with probability approaching one, if $L$ is large enough.
 
We have 
\begin{eqnarray*}
&& \sum_{i=1}^m\sum_{t=1}^T (y_{it}-\vtheta_i\trans\vB(\vX_{it}\trans\vbeta_i))^2- \sum_{i=1}^m\sum_{t=1}^T (y_{it}-\vtheta_{0i}\trans\vB(\vX_{it}\trans\vbeta_{0i}))^2\nonumber\\
&=&\sum_{i,t}( \vtheta_i\trans\vB(\vX_{it}\trans\vbeta_i)-\vtheta_{0i}\trans\vB(\vX_{it}\trans\vbeta_{0i}))^2-2(\epsilon_{it}-r_{it})(\vtheta_i\trans\vB(\vX_{it}\trans\vbeta_i)-\vtheta_{0i}\trans\vB(\vX_{it}\trans\vbeta_{0i})),\nonumber\\
\end{eqnarray*}
where $r_{it}=\vtheta_{0i}\trans\vB(\vX_{it}\trans\vbeta_{0i})-g(\vX_{it}\trans\vbeta_{0i})$ with $|r_{it}|\le  CK^{-2}$.

Furthermore,
\begin{eqnarray*}
&&\sum_{i,t}(\vtheta_i\trans\vB(\vX_{it}\trans\vbeta_i)-\vtheta_{0i}\trans\vB(\vX_{it}\trans\vbeta_{0i}))^2\\
&=&\sum_{i,t} \left( (\vtheta_i-\vtheta_{0i})\trans\vB(\vX_{it}\trans\vbeta_i)+\vtheta_{0i}\trans\vB'(\vX_{it}\vbeta_i^*)\wb\vX_{it}\trans(\wb\vbeta_i-\wb\vbeta_{0i}) \right)^2\\
&=&T (\vtheta_1\trans-\vtheta_{01}\trans,\wb\vbeta_1\trans-\wb\vbeta_{01}\trans,\ldots,\vtheta_m\trans-\vtheta_{0m}\trans,\wb\vbeta_m\trans-\wb\vbeta_{0m}\trans)\cdot\\
&&\left(\begin{array}{cccc}
    \wt\vA_{11}&\vnull&\cdots&\vnull\\
    \vnull&  \wt\vA_{22}&\cdots&\vnull\\
    \vdots&\vdots&\vdots&\vdots\\
    \vnull&\vnull&\cdots&\wt\vA_{mm}
	\end{array}\right) \cdot \left(\begin{array}{c}
	\vtheta_1-\vtheta_{01}\\
	\wb\vbeta_1-\wb\vbeta_{01}\\
	\vdots\\
	\vtheta_m-\vtheta_{m1}\\	
	\wb\vbeta_m-\wb\vbeta_{m1}
	\end{array}\right)\\
&=&T ((\vxi\trans-\vxi_0\trans) \vD^{G_1}, (\veta\trans-\veta_0\trans)\vD^{G_2})\\
&&   \left(\begin{array}{cc}
      (\vD^{G_1})^{-1}&\vnull\\
      \vnull & (\vD^{G_2})^{-1}
         \end{array}\right)
   \left(\begin{array}{ccccc}
   		(\vJ_1^{G_1})\trans &\vnull &\cdots &(\vJ_m^{G_1})\trans &\vnull\\
   		\vnull & (\vJ_1^{G_2})\trans &\cdots &\vnull &(\vJ_m^{G_2})\trans 
   		\end{array}\right)\\
&&\left(\begin{array}{cccc}
     \wt\vA_{11}&\vnull&\cdots&\vnull\\
    \vnull&  \wt\vA_{22}&\cdots&\vnull\\
    \vdots&\vdots&\vdots&\vdots\\
    \vnull&\vnull&\cdots&\wt\vA_{mm}
	\end{array}\right)\cdot 
   \left(\begin{array}{cc}
        \vJ_1^{G_1}&\vnull\\
        \vnull & \vJ_1^{G_2}\\
        \vdots &\vdots\\
        \vJ_m^{G_1}&\vnull\\
        \vnull & \vJ_m^{G_2}\\
   		\end{array}\right)
	   \left(\begin{array}{cc}
      (\vD^{G_1})^{-1}&\vnull\\
      \vnull & (\vD^{G_2})^{-1}
         \end{array}\right)
       \left(\begin{array}{c}
           \vD^{G_1}(\vxi-\vxi_0)\\
           \vD^{G_2}(\veta-\veta_0)
           \end{array}\right),
\end{eqnarray*}
where 
$$\wt\vA_{ii'}=\frac{1}{T}\sum_{t=1}^T\left[\left(\begin{array}{c}
	\vB(\vX_{it}\trans\vbeta_{i})\\
	\vtheta_{0i}\trans \vB'(\vX_{it}\trans\vbeta_{i}^*)\wb\vX_{it}
	\end{array}\right)
	\left(\begin{array}{cc}
     	\vB\trans(\vX_{it}\trans\vbeta_i) & \vtheta_{0i}\trans\vB'(\vX_{it}\trans\vbeta_{i}^*)\wb\vX_{it}\trans
	\end{array}\right)\right]
	, 1\le i,i'\le m,$$
$\vB'(.)=(B_1'(.),\ldots,B_K'(.))\trans$ are the first derivatives of the basis functions and $\vbeta_i^*$ lies between $\vbeta_{0i}$ and $\vbeta_i$.

By Lemma \ref{lem:eigenAtilde}, eigenvalues of $\wt\vA_{ii}$ are bounded and bounded away from zero, with probability approaching one. Furthermore, it is easy to directly verify that 
\begin{equation}\label{eqn:vo}
  \vO:= \left(\begin{array}{cc}
        \vJ_1^{G_1}&\vnull\\
        \vnull & \vJ_1^{G_2}\\
        \vdots &\vdots\\
        \vJ_m^{G_1}&\vnull\\
        \vnull & \vJ_m^{G_2}\\
   		\end{array}\right)
	   \left(\begin{array}{cc}
      (\vD^{G_1})^{-1}&\vnull\\
      \vnull & (\vD^{G_2})^{-1}
         \end{array}\right)
\end{equation}
is an orthonormal matrix (that is, $\vO\trans\vO=\vI$). Thus 
\begin{equation}\label{eqn:comb1}
\sum_{i,t}(\vtheta_i\trans\vB(\vX_{it}\trans\vbeta_i)-\vtheta_{0i}\trans\vB(\vX_{it}\trans\vbeta_{0i}))^2\asymp T(\|\vD^{G_1}(\vxi-\vxi_0)\|^2+\|\vD^{G_2}(\veta-\veta_0)\|^2)=T(\|\vtheta-\vtheta_0\|^2+\|\vbeta-\vbeta_0\|^2).
\end{equation}

Now consider the term 
$(\epsilon_{it}-r_{it})(\vtheta_i\trans\vB(\vX_{it}\trans\vbeta_i)-\vtheta_{0i}\trans\vB(\vX_{it}\trans\vbeta_{0i}))$.
We have
\begin{eqnarray*}
&&\sum_{i,t}\epsilon_{it}(\vtheta_i\trans\vB(\vX_{it}\trans\vbeta_i)-\vtheta_{0i}\trans\vB(\vX_{it}\trans\vbeta_{0i}))\\
&=&\sum_{t}(\vtheta_1\trans-\vtheta_{01}\trans,\wb\vbeta_1\trans-\wb\vbeta_{01}\trans,\ldots,\vtheta_m\trans-\vtheta_{0m}\trans,\wb\vbeta_m\trans-\wb\vbeta_{0m}\trans)\cdot\\
&&\left(\begin{array}{c}
	\vB(\vX_{1t}\trans\vbeta_1)\epsilon_{1t}\\
	\vtheta_{01}\trans\vB'(\vX_{1t}\trans\vbeta_1^*)\wb\vX_{1t}\epsilon_{1t}\\
	\vdots\\
	\vB(\vX_{mt}\trans\vbeta_m)\epsilon_{mt}\\
	\vtheta_{0m}\trans\vB'(\vX_{mt}\trans\vbeta_m^*)\wb\vX_{mt}\epsilon_{mt}
	\end{array}\right)\\
&=&\sum_t ((\vxi\trans-\vxi_0\trans) \vD^{G_1}, (\veta\trans-\veta_0\trans)\vD^{G_2})
  \vO\trans\left(\begin{array}{c}
	\vB(\vX_{1t}\trans\vbeta_1)\epsilon_{1t}\\
	\vtheta_{01}\trans\vB'(\vX_{1t}\trans\vbeta_1^*)\wb\vX_{1t}\epsilon_{1t}\\
	\vdots\\
	\vB(\vX_{mt}\trans\vbeta_m)\epsilon_{mt}\\
	\vtheta_{0m}\trans\vB'(\vX_{mt}\trans\vbeta_m^*)\wb\vX_{mt}\epsilon_{mt}
	\end{array}\right)\\
&\le& \sqrt{\|\vD^{G_1}(\vxi-\vxi_0)\|^2+\|\vD^{G_2}(\veta-\veta_0)\|^2} \left\|\sum_t\vO\trans \left(\begin{array}{c}
	\vB(\vX_{1t}\trans\vbeta_1)\epsilon_{1t}\\
	\vtheta_{01}\trans\vB'(\vX_{1t}\trans\vbeta_1^*)\wb\vX_{1t}\epsilon_{1t}\\
	\vdots\\
	\vB(\vX_{mt}\trans\vbeta_m)\epsilon_{mt}\\
	\vtheta_{0m}\trans\vB'(\vX_{mt}\trans\vbeta_m^*)\wb\vX_{mt}\epsilon_{mt}
	\end{array}\right)\right\|,
\end{eqnarray*}
where $\vO$ is as defined in \eqref{eqn:vo}. We  have
\begin{eqnarray}\label{eqn:traceinq}
&&E\left\|\sum_t\vO\trans \left(\begin{array}{c}
	\vB(\vX_{1t}\trans\vbeta_1)\epsilon_{1t}\\
	\vtheta_{01}\trans\vB'(\vX_{1t}\trans\vbeta_1^*)\wb\vX_{1t}\epsilon_{1t}\\
	\vdots\\
	\vB(\vX_{mt}\trans\vbeta_m)\epsilon_{mt}\\
	\vtheta_{0m}\trans\vB'(\vX_{mt}\trans\vbeta_m^*)\wb\vX_{mt}\epsilon_{mt}
	\end{array}\right)\right\|^2\nonumber\\
&=&{\rm tr}\left( \sum_{1\le t,t'\le T}\vO\trans\cdot \left[\begin{array}{ccc}
		\vA_{11,|t-t'|}\sigma_{11,|t-t'|}&\cdots &\vA_{1m,|t-t'|}\sigma_{1m,|t-t'|}\\
		\vdots&\vdots&\vdots\\
		\vA_{m1,|t-t'|}\sigma_{m1,|t-t'|}&\cdots &\vA_{mm,|t-t'|}\sigma_{mm,|t-t'|}\end{array}\right]\vO\right)\nonumber\\
&=&{\rm tr}\left( \sum_{1\le t,t'\le T} \left[\begin{array}{ccc}
		\vA_{11,|t-t'|}\sigma_{11,|t-t'|}&\cdots &\vA_{1m,|t-t'|}\sigma_{1m,|t-t'|}\\
		\vdots&\vdots&\vdots\\
		\vA_{m1,|t-t'|}\sigma_{m1,|t-t'|}&\cdots &\vA_{mm,|t-t'|}\sigma_{mm,|t-t'|}\end{array}\right]\vO\vO\trans\right)\nonumber\\
&\le&{\rm tr}(\vO\vO\trans)\cdot \left\| \sum_{1\le t,t'\le T} \left[\begin{array}{ccc}
		\vA_{11,|t-t'|}\sigma_{11,|t-t'|}&\cdots &\vA_{1m,|t-t'|}\sigma_{1m,|t-t'|}\\
		\vdots&\vdots&\vdots\\
		\vA_{m1,|t-t'|}\sigma_{m1,|t-t'|}&\cdots &\vA_{mm,|t-t'|}\sigma_{mm,|t-t'|}\end{array}\right] \right\|_{op},
\end{eqnarray}
where $\sigma_{ii',|t-t'|}=Cov(\epsilon_{it},\epsilon_{i't'})$,
$$\vA_{ii',|t-t'|}=E\left[\left(\begin{array}{c}
	\vB(\vX_{it}\trans\vbeta_{i})\\
	\vtheta_{0i}\trans \vB'(\vX_{it}\trans\vbeta_{i}^*)\wb\vX_{it}
	\end{array}\right)
	\left(\begin{array}{cc}
     	\vB\trans(\vX_{i't'}\trans\vbeta_{i'}) & \vtheta_{0i'}\trans\vB'(\vX_{i't'}\trans\vbeta_{i'}^*)\wb\vX_{i't'}\trans
	\end{array}\right)\right],
$$
and the last step above uses von Neumann's trace inequality \citep{mirsky75}. By Lemma \ref{lem:eigenAtilde2} and that ${\rm tr}(\vO\vO\trans)=H_1+H_2$ (note $\vO\trans\vO=\vI_{H_1+H_2}$), we have
\begin{eqnarray}\label{eqn:comb2}
&&\sum_{i,t}\epsilon_{it}(\vtheta_i\trans\vB(\vX_{it}\trans\vbeta_i)-\vtheta_{0i}\trans\vB(\vX_{it}\trans\vbeta_{0i}))\nonumber\\
&=&O_p(\sqrt{(\|\vtheta-\vtheta\|^2+\|\vbeta-\vbeta_0\|^2)(H_1+H_2)T}).
\end{eqnarray}
Finally, using Cauchy-Schwarz inequality
\begin{eqnarray}\label{eqn:comb3}
&&\sum_{i,t}r_{it}(\vtheta_i\trans\vB(\vX_{it}\trans\vbeta_i)-\vtheta_{0i}\trans\vB(\vX_{it}\trans\vbeta_{0i}))\nonumber\\
&=&C\sqrt{mT}K^{-2}\cdot O_p(\sqrt{T(\|\vtheta-\vtheta\|^2+\|\vbeta-\vbeta_0\|^2)} )
\end{eqnarray}
Combining \eqref{eqn:comb1}--\eqref{eqn:comb3},
\begin{eqnarray*}
&& \sum_{i=1}^m\sum_{t=1}^T (y_{it}-\vtheta_i\trans\vB(\vX_{it}\trans\vbeta_i))^2- \sum_{i=1}^m\sum_{t=1}^T (y_{it}-\vtheta_{0i}\trans\vB(\vX_{it}\trans\vbeta_{0i}))^2>0  
\end{eqnarray*}
with probability approaching one, if $\|\vbeta-\vbeta_0\|^2+\|\vtheta-\vtheta_0\|^2= Lr_T^2$ with $L$ sufficiently large. Thus there is a local minimizer $(\wh\vtheta,\wh\vbeta)$ with 
$\|\wh\vbeta-\vbeta_0\|+\|\wh\vtheta-\vtheta_0\|=O_p(r_T)$.

\bigskip

\noindent\textit{STEP 2. Proof of convergence rate of $\wh\vbeta$ and its asymptotic normality}.

Let $\vPi_i$ be $T\times K$ matrices,  $i=1,\ldots,m$, with rows $\vPi_{it}\trans=\vB(\vX_{it}\trans\vbeta_{0i})$. Define $\vV_{it}=g_{0i}'(\vX_{it}\trans\vbeta_{0i})\wb\vX_{it}$, $\vP_i=\vPi_i(\vPi_i\trans\vPi_i)^{-1}\vPi_i\trans$ with rows $\vP_{it}\trans=\vPi_{it}\trans(\vPi_i\trans\vPi_i)^{-1}\vPi_i\trans$. We write, for any $(\vtheta,\vbeta)$ with $\|\vtheta-\vtheta_0\|^2\le  Cr_T^2$ and $\|\vbeta-\vbeta_0\|^2\le CH_2/T$,
\begin{eqnarray*}
&&\sum_{i,t}(y_{it}-\vtheta_i\trans\vB(\vX_{it}\trans\vbeta_i))^2\\
&=&\sum_{i,t}(\epsilon_{it}+g_{0i}(\vX_{it}\trans\vbeta_i)-\vtheta_i\trans\vB(\vX_{it}\trans\vbeta_i))^2\\
&=&\sum_{i,t}(\epsilon_{it}-\vPi_{it}\trans(\vtheta_i-\vtheta_{0i})-\vV_{it}\trans(\wb\vbeta_i-\wb\vbeta_{0i})-R_{it})^2,
\end{eqnarray*}
where 
\begin{eqnarray*}
&&R_{it}\\
&=&\vtheta_i\trans\vB(\vX_{it}\trans\vbeta_i)-g_{0i}(\vX_{it}\trans\vbeta_{0i})-(\vtheta_i-\vtheta_{0i})\trans\vB(\vX_{it}\trans\vbeta_{0i})-g_{0i}'(\vX_{it}\trans\vbeta_{0i})\wb\vX_{it}\trans(\wb\vbeta_i-\wb\vbeta_{0i})\\
&=&\vtheta_{0i}\trans\vB(\vX_{it}\trans\vbeta_{0i})-g_{0i}(\vX_{it}\trans\vbeta_{0i})+\vtheta_i\trans(\vB(\vX_{it}\trans\vbeta_i)-\vB(\vX_{it}\trans\vbeta_{0i}))-g_{0i}'(\vX_{it}\trans\vbeta_{0i})\wb\vX_{it}\trans(\wb\vbeta_i-\wb\vbeta_{0i})\\
&=&\vtheta_{0i}\trans\vB(\vX_{it}\trans\vbeta_{0i})-g_{0i}(\vX_{it}\trans\vbeta_{0i})+(\vtheta_i-\vtheta_{0i})\trans(\vB(\vX_{it}\trans\vbeta_{i})-\vB(\vX_{it}\trans\vbeta_{0i}))+\vtheta_{0i}\trans(\vB(\vX_{it}\trans\vbeta_{i})-\vB(\vX_{it}\trans\vbeta_{0i}))\\
&&-g_{0i}'(\vX_{it}\trans\vbeta_{0i})\wb\vX_{it}\trans(\wb\vbeta_i-\wb\vbeta_{0i})\\
&=&\left\{\vtheta_{0i}\trans\vB(\vX_{it}\trans\vbeta_{0i})-g_{0i}(\vX_{it}\trans\vbeta_{0i})\right\}+\left\{(\vtheta_i-\vtheta_{0i})\trans(\vB(\vX_{it}\trans\vbeta_{i})-\vB(\vX_{it}\trans\vbeta_{0i}))\right\}\\
&&+\left\{\vtheta_{0i}\trans(\vB(\vX_{it}\trans\vbeta_{i})-\vB(\vX_{it}\trans\vbeta_{0i})-\vB'(\vX_{it}\trans\vbeta_{0i})\wb\vX_{it}\trans(\wb\vbeta_i-\wb\vbeta_{0i})\right\}\\
&&+\left\{(\vtheta_{0i}\trans\vB'(\vX_{it}\trans\vbeta_{0i})-g_{0i}'(\vX_{it}\trans\vbeta_{0i}))\wb\vX_{it}\trans(\wb\vbeta_i-\wb\vbeta_{0i})\right\}\\
&=&R_{it1}+R_{it2}(\vtheta_i,\vbeta_i),
\end{eqnarray*}
where $R_{it1}=\vtheta_0\trans\vB(\vX_{it}\trans\vbeta_{0i})-g_{0i}(\vX_{it}\trans\vbeta_{0i})$ and $R_{it2}(\vtheta_i,\vbeta_i)$ contains all other terms above. It is easy to see $R_{it2}(\vtheta_i,\vbeta_{0i})=0$.  In the decomposition above $R_{it2}$ consists of three terms, which we denote by $R_{it2,1}$, $R_{it2,2}$ and $R_{it2,3}$, respectively (omitting the dependence in $\vtheta,\vbeta$ for simplicity of notation). 
Using $\|\vtheta-\vtheta_0\|^2+\|\vbeta-\vbeta_0\|^2\le  Cr_T^2$, we can easily show 
\begin{eqnarray*}
\sum_{i,t}R_{it2,1}^2&=& O_p(Tr_T^4K^3),\\
\sum_{i,t}R_{it2,2}^2&=& O_p(Tr_T^4p),\\
\sum_{i,t}R_{it2,3}^2&=& O_p(Tr_T^2K^{-2}),
\end{eqnarray*}
and  thus 
\be\label{eqn:Rit2}
\sum_{i,t}R_{it2}^2=O_p\left(Tr_T^4(K^3+p)+Tr_T^2K^{-2}\right).
\ee
We then orthogonalize the parametric part with respect to the nonparametric part by  writing 
\begin{eqnarray*}
&&\sum_{i,t}(\epsilon_{it}-\vPi_{it}\trans(\vtheta_i-\vtheta_{0i})-\vV_{it}\trans(\wb\vbeta_i-\wb\vbeta_{0i})-R_{it})^2\\
&=&\sum_{i,t}(\epsilon_{it}-\vPi_{it}\trans(\valpha_i-\valpha_{0i})-(\vV_{it}-\vV_i\trans\vP_{it})\trans(\wb\vbeta_i-\wb\vbeta_{0i})-R_{it1}-R_{it2}(\calM_i(\valpha_i,\vbeta_i))^2,\\
\end{eqnarray*}
where $\valpha_i=\vtheta_i+(\vPi_i\trans\vPi_i)^{-1}\vPi_i\trans\vV_i\wb\vbeta_i$, $\valpha_{0i}=\vtheta_{0i}+(\vPi_i\trans\vPi_i)^{-1}\vPi_i\trans\vV_i\wb\vbeta_{0i}$, $\vV_i=(\vV_{it},\ldots,\vV_{iT})\trans$, and $\calM_i$ is the one-to-one mapping that maps $(\valpha_i,\vbeta_i)$ to $(\vtheta_i,\vbeta_i)$. Below we write $R_{it2}(\calM_i(\valpha_i,\vbeta_i))$ as $R_{it2}$, $R_{it2}(\calM_i(\wh\valpha_i,\wh\vbeta_i))$ as $\wh R_{it2}$,  and note $R_{it2}(\calM_i(\wh\valpha_i,\vbeta_0))=0$.
Then,
\begin{eqnarray}\label{eqn:diff0}
0
&\ge&\sum_{i,t}(\epsilon_{it}-\vPi_{it}\trans\wh\valpha_i-(\vV_{it}-\vV_i\trans\vP_{it})\trans(\wh{\wb\vbeta}_i-\wb\vbeta_{0i})-R_{it1}-R_{it2})^2\nonumber\\
&&-\sum_{i,t}(\epsilon_{it}-\vPi_{it}\trans\wh\valpha_i-R_{it1})^2\nonumber\\
&=&\sum_{i,t}(\wh\veta-\veta_0)\trans(\vJ_i^{G_2})\trans(\vV_{it}-\vV_i\trans\vP_{it})(\vV_{it}\trans-\vP_{it}\trans\vV_i)\vJ_i^{G_2}(\wh\veta-\veta_0)+\sum_{i,t}\wh R_{it2}^2\nonumber\\
&&-2\sum_{i,t}\left((\wh\veta-\veta_0)\trans(\vJ_i^{G_2})\trans(\vV_{it}-\vV_i\trans\vP_{it})+\wh R_{it2}\right)(\epsilon_{it}-\vPi_{it}\trans\valpha_i-R_{it1})\nonumber\\
&&+2\sum_{i,t}\wh R_{it2}\cdot(\wh\veta-\veta_0)\trans(\vJ_i^{G_2})\trans(\vV_{it}-\vV_i\trans\vP_{it})\nonumber\\
&=&\sum_{i,t}(\wh\veta-\veta_0)\trans(\vJ_i^{G_2})\trans(\vV_{it}-\vV_i\trans\vP_{it})(\vV_{it}\trans-\vP_{it}\trans\vV_i)\vJ_i^{G_2}(\wh\veta-\veta_0)\nonumber\\
&&-2\sum_{i,t}(\wh\veta-\veta_0)\trans(\vJ_i^{G_2})\trans(\vV_{it}-\vV_i\trans\vP_{it})\epsilon_{it}\nonumber\\
&&-2\sum_{i,t}\wh R_{it2}\epsilon_{it}\nonumber\\
&&+\sum_{i,t}\wh R_{it2}^2+2\sum_{i,t}\left((\wh\veta-\veta_0)\trans(\vJ_i^{G_2})\trans(\vV_{it}-\vV_i\trans\vP_{it})+\wh R_{it2}\right)(\vPi_{it}\trans\wh\valpha_i+R_{it1})\nonumber\\
&&+2\sum_{i,t}\wh R_{it2}\cdot(\wh\veta-\veta_0)\trans(\vJ_i^{G_2})\trans(\vV_{it}-\vV_i\trans\vP_{it}).
\end{eqnarray}
The first term above is
\bse
&&\sum_{i,t}(\wh\veta-\veta_0)(\vJ_i^{G_2})\trans(\vV_{it}-\vV_i\trans\vP_{it})(\vV_{it}\trans-\vP_{it}\trans\vV_i)\vJ_i^{G_2}(\wh\veta-\veta_0)\\
&=&T(\wh\veta-\veta_0)\trans\vD^{G_2}\vO_2\trans\left(\begin{array}{cccc}
					\wh\vC_{11}&\vnull&\cdots&\vnull\\
					\vnull&\wh\vC_{22}&\cdots&\vnull\\					
					\vdots&\vdots&\vdots&\vdots\\
					\vnull&\vnull&\cdots&\wh\vC_{mm}
					\end{array}\right)\vO_2\vD^{G_2}(\wh\veta-\veta_0),					
\ese
where $\vO_2=\left(\begin{array}{c} \vJ_1^{G_2}\\\vdots\\\vJ_m^{G_2}\end{array}\right)(\vD^{G_2})^{-1}$ is an $mp\times H_2$ orthonormal matrix, and $\wh\vC_{ii}=\sum_{t}(\vV_{it}-\vV_i\trans\vP_{it})(\vV_{it}-\vV_i\trans\vP_{it})\trans/T$. 

Let $\vC_{ii}=E[(g_{0i}'(\vX_{it}\trans\vbeta_{0i}))^2(\wb\vX_{it}-E[\wb\vX_{it}|\vX_{it}\vbeta_{0i}])^{\otimes 2}]$. Lemma \ref{lem:eigenC} shows that $\max_i\|\wh\vC_{ii}-\vC_{ii}\|_{op}=o_p(1)$. 
Based on this, we have the first term in (\ref{eqn:diff0})  is bounded below by $CT\|\vD^{G_2}(\wh\veta-\veta_0)\|^2$.

Now consider the second term in (\ref{eqn:diff0}). We have
\bse
&&\sum_{i,t}(\wh\veta-\veta_0)(\vJ_i^{G_2})\trans(\vV_{it}-\vV_i\trans\vP_{it})\epsilon_{it}\\
&\le&\|\vD^{G_2}(\wh\veta-\veta_0)\|\cdot \left\| \sum_t\vO_2\trans	\left(\begin{array}{c}
						(\vV_{1t}-\vV_1\trans\vP_{1t})\epsilon_{1t}\\			
						\vdots\\
						(\vV_{mt}-\vV_m\trans\vP_{mt})\epsilon_{mt}
						\end{array}\right)\right\|.
\ese
We write 
\bse
&& \sum_t\left(\begin{array}{c}
						(\vV_{1t}-\vV_1\trans\vP_{1t})\epsilon_{1t}\\			
						\vdots\\
						(\vV_{mt}-\vV_m\trans\vP_{mt})\epsilon_{mt}
						\end{array}\right)\\
&=& \sum_t\left(\begin{array}{c}
						(\vV_{1t}-\vPhi_{1t})\epsilon_{1t}\\			
						\vdots\\
						(\vV_{mt}-\vPhi_{mt})\epsilon_{mt}
						\end{array}\right)+\left(\begin{array}{c}
		((\vI-\vP_1)\vPhi_1-\vP_1(\vV_1-\vPhi_1))\trans\vepsilon_1\\
		\vdots\\
		((\vI-\vP_m)\vPhi_m-\vP_m(\vV_m-\vPhi_m))\trans\vepsilon_m
			\end{array}\right),
\ese
where $\vepsilon_i=(\epsilon_{i1},\ldots,\epsilon_{iT})\trans$. 

The covariance matrix of $	\left(\begin{array}{c}
						(\vV_{1t}-\vPhi_{1t})\epsilon_{1t}\\			
						\vdots\\
						(\vV_{mt}-\vPhi_{mt})\epsilon_{mt}
						\end{array}\right)$ is given by
$$\sum_{1\le t,t'\le T}\left[\begin{array}{ccc}
		\vC_{11,|t-t'|}\sigma_{11,|t-t'|}&\cdots &\vC_{1m,|t-t'|}\sigma_{1m,|t-t'|}\\
		\vdots&\vdots&\vdots\\
		\vC_{m1,|t-t'|}\sigma_{m1,|t-t'|}&\cdots &\vC_{mm,|t-t'|}\sigma_{mm,|t-t'|}\end{array}\right],$$
with $\vC_{ii',|t-t'|}=E[g_{0i}'(\vX_{it}\trans\vbeta_{0i})g_{0i'}'(\vX_{i't'}\trans\vbeta_{0i'})(\wb\vX_{it}-E[\wb\vX_{it}|\vX_{it}\trans\vbeta_{0i}])(\wb\vX_{i't'}-E[\wb\vX_{i't'}|\vX_{i't'}\trans\vbeta_{0i'}])\trans]$.
Using the geometric mixing rate, and similar to the proof of Lemma \ref{lem:eigenAtilde2}, it can be shown that the matrix above has eigenvalues of order $O_p(T)$.

Furthermore, we can bound the  largest eigenvalue of 
$$E\left[\left(\begin{array}{c}
		((\vI-\vP_1)\vPhi_1-\vP_1(\vV_1-\vPhi_1))\trans\vepsilon_1\\
		\vdots\\
		((\vI-\vP_m)\vPhi_m-\vP_m(\vV_m-\vPhi_m))\trans\vepsilon_m
			\end{array}\right)^{\otimes 2}\right].
$$
Denoting $\vE_i=(\vI-\vP_i)\vPhi_i+\vP_i(\vV_i-\vPhi_i)$, in Lemma \ref{lem:eigenC} we have shown that $\max_i\|\vE_i\|^2=O_p(TK^{-4}+Kp\log T)$. We have 
\bse
&&E\left[\left.\left(\begin{array}{c}
		\vE_1\trans\vepsilon_1\\
		\vdots\\
		\vE_m\trans\vepsilon_m
			\end{array}\right)^{\otimes 2}\right|\{\vX_{it}\}\right]\\
&=&\left(\begin{array}{ccc}
         \vE_1\trans E[\vepsilon_1\vepsilon_1\trans]\vE_1 &\ldots& \vE_1\trans E[\vepsilon_1\vepsilon_m\trans]\vE_m\\
         \vdots&\vdots&\vdots\\
         \vE_m\trans E[\vepsilon_m\vepsilon_1\trans]\vE_1 &\ldots& \vE_m\trans E[\vepsilon_m\vepsilon_m\trans]\vE_m\\
        \end{array}\right).
\ese
Note $\|\vE_i\trans E[\vepsilon_i\vepsilon_{i'}\trans]\vE_{i'}\|_{op}\le \|E[\vepsilon_i\vepsilon_{i'}\trans]\|_{op} \|\vE_i\|\|\vE_{i'}\|$. Since for any $\vv,\vu\in R^T$, $\vu\trans E[\vepsilon_i\vepsilon_i\trans]\vv\le (\vu\trans E[\vepsilon_i\vepsilon_i\trans]\vu+\vv\trans E[\vepsilon_{i'}\vepsilon_{i'}\trans]\vv)/2$, and (by assumption (C2)) $\|E[\vepsilon_{i}\vepsilon_{i}\trans]\|_{op}\le M$, we have $\|E[\vepsilon_{i}\vepsilon_{i'}\trans]\|_{op}\le M$ for all $(i,i')$. Furthermore, $\max_{i,i'}\|\vE_i\|\|\vE_{i'}\|=O_p(TK^{-4}+Kp\log T)$. Thus $\|\vE_i\trans E[\vepsilon_i\vepsilon_{i'}\trans]\vE_{i'}\|_{op}=O_p(TK^{-4}+Kp\log T)$, uniformly over $(i,i')$. Now for $\vv=(\vv_1\trans,\ldots,\vv_m\trans)\trans\in R^{mp}$, 
\bse
&&\vv\trans\left(\begin{array}{ccc}
         \vE_1\trans E[\vepsilon_1\vepsilon_1\trans]\vE_1 &\ldots& \vE_1\trans E[\vepsilon_1\vepsilon_m\trans]\vE_m\\
         \vdots&\vdots&\vdots\\
         \vE_m\trans E[\vepsilon_m\vepsilon_1\trans]\vE_1 &\ldots& \vE_m\trans E[\vepsilon_m\vepsilon_m\trans]\vE_m\\
        \end{array}\right)\vv\\
&=&\sum_{i,i'}\vv_i\trans \vE_i\trans E[\vepsilon_i\vepsilon_{i'}\trans]\vE_{i'}\vv_{i'}\\
&\le& \sum_{i,i'}\|\vv_i\|\|\vE_i\trans E[\vepsilon_i\vepsilon_{i'}\trans]\vE_{i'}\|_{op}\|\vv_{i'}\|,
\ese
which is bounded by the largest eigenvalue of the $m\times m$ matrix with entries $\lambda_{\max}(\vE_i\trans E[\vepsilon_i\vepsilon_{i'}\trans]\vE_{i'}), i,i'=1,\ldots,m$. This matrix has eigenvalues bounded by $O_p(mTK^{-4}+mKp\log T)=o_p(T)$ by the Gershgorin circle theorem.

Using the trace inequality as in (\ref{eqn:traceinq}), we get 
\be\label{eqn:second}
&&E\left[ \left\| \sum_t\vO_2\trans	\left(\begin{array}{c}
						(\vV_{1t}-\vV_1\trans\vP_{1t})\epsilon_{1t}\\			
						\vdots\\
						(\vV_{mt}-\vV_m\trans\vP_{mt})\epsilon_{mt}
						\end{array}\right)\right\|^2\right]=O(H_2T)
\ee
and thus the second term in (\ref{eqn:diff0}) is $O_p(\sqrt{H_2T}\|\vD^{G_2}(\wh\vbeta-\vbeta_0)\|)$. For the rest of the terms in (\ref{eqn:diff0}),
we have, using (\ref{eqn:Rit2}),
\begin{eqnarray*}
&&(\sum_{i,t}\wh R_{it2}\epsilon_{it})^2=O_p(Tr_T^4(K^3+p)+Tr_T^2K^{-2}),\\
&&\sum_{i,t}(\wh\veta-\veta_0)\trans(\vJ_i^{G_2})\trans(\vV_{it}-\vV_i\trans\vP_{it})\vPi_{it}\trans\wh\valpha_i=\sum_{i}(\wh\veta-\veta_0)\trans(\vJ_i^{G_2})\trans(\vV_{i}-\vP_i\vV_i)\vPi_{it}\trans\wh\valpha_i=0,\\
&&(\sum_{i,t}\wh R_{it2}R_{it1})^2\le (\sum_{i,t}\wh R_{it2}^2)(\sum_{i,t}R_{it1}^2)=O_p((Tr_T^4(K^3+p)+Tr_T^2K^{-2})mTK^{-4}),\\
&&(\sum_{i,t}\wh R_{it2}\vPi_{it}\trans\wh\valpha_i)^2\le (\sum_{i,t}\wh R_{it2}^2)(\sum_{i,t}(\vPi_{it}\trans\wh\valpha_i)^2)=O_p((Tr_T^4(K^3+p)+Tr_T^2K^{-2})Tr_T^2),\\
&&(\sum_{i,t}\wh R_{it2}(\wh\veta-\veta_0)\trans\vJ_i^{G_2}(\vV_{it}-\vV_i\trans\vP_{it}))^2\\
&&\le (\sum_{i,t}\wh R_{it2}^2)\left\|(\vD^{G_2}(\wh\veta-\veta_0))\trans\vO_2\trans\left(\begin{array}{c}
			\vV_1\trans-\vV_1\trans\vP_1\trans\\
			\vdots\\
			\vV_m\trans-\vV_m\trans\vP_m\trans \end{array}\right)\right\|^2\\	&&=O_p((Tr_T^4(K^3+p)+Tr_T^2K^{-2})mpH_2).
\end{eqnarray*}
All these terms are order $o_p(1)$ by our assumptions.
Finally, consider the term 
\begin{eqnarray*}
&&(\sum_{i,t}R_{it1}(\veta-\veta_0)\trans(\vJ_i^{G_2})\trans(\vV_{it}-\vV_i\trans
\vP_{it}))^2\\
&=&O_p(H_2/T)\|\sum_i\vR_{i1}\trans (\vV_i-\vP_i\vV_i)\|^2,
\end{eqnarray*} 
where $\vR_{i1}=(R_{i11},\ldots,R_{iT1})\trans$.
Again, using $\vV_i=(\vV_i-\vPhi_i)+\vPhi_i$,
\bse
&&(\sum_{i,t}R_{it1}(\veta-\veta_0)\trans(\vJ_i^{G_2})\trans(\vV_{it}-\vV_i\trans
\vP_{it}))^2\\
&=&O_p(H_2/T)\left(\|\sum_i\vR_{i1}\trans(\vI-\vP_i)\vPhi_i\|^2+\|\sum_i\vR_{i1}\vP_i(\vV_i-\vPhi_i)\|^2+\|\sum_i\vR_{i1}(\vV_i-\vPhi_i)\|^2\right)\\
&=&O_p(H_2/T)\left(O_p(mTK^{-4}\cdot mTK^{-4})+O_p(mTK^{-4}\cdot mKp\log T)+O_p(m^2pTK^{-4})\right)\\
&=&o_p(1).
\ese
Summarizing the bounds for different terms in (\ref{eqn:diff0}), we get
\bse
\|\vD^{G_2}(\wh\veta-\veta_0)\|^2+\|\vD^{G_2}(\wh\veta-\veta_0)\|O_p(\sqrt{H_2/T})+o_p(1/T)\le 0
\ese
Completing the squares, we get 
\bse
(\|\vD^{G_2}(\wh\veta-\veta_0)\|+O_p(\sqrt{H_2/T}))^2=O_p(H_2/T)
\ese
which in turn implies $\|\wh\vbeta-\vbeta_0\|=\|\vD^{G_2}(\wh\veta-\veta_0)\|=O_p(\sqrt{H_2/T})$.

To get asymptotic normality, we similarly write
\begin{eqnarray}\label{eqn:resid}
&&\sum_{i,t}(\epsilon_{it}-\vPi_{it}\trans\valpha_i-(\vV_{it}-\vV_i\trans\vP_{it})\trans(\wb\vbeta_i-\wb\vbeta_{0i})-R_{it1}-R_{it2})^2\nonumber\\
&&-\sum_{i,t}(\epsilon_{it}-\vPi_{it}\trans\valpha_i-R_{it1})^2\nonumber\\
&=&\sum_{i,t}(\veta-\veta_0)\trans(\vJ_i^{G_2})\trans(\vV_{it}-\vV_i\trans\vP_{it})(\vV_{it}\trans-\vP_{it}\trans\vV_i)\vJ_i^{G_2}(\veta-\veta_0)\nonumber\\
&&-2\sum_{i,t}(\veta-\veta_0)\trans(\vJ_i^{G_2})\trans(\vV_{it}-\vV_i\trans\vP_{it})\epsilon_{it}\nonumber\\
&&-2\sum_{i,t}R_{it2}\epsilon_{it}\nonumber\\
&&+\sum_{i,t}R_{it2}^2+2\sum_{i,t}\left((\veta-\veta_0)\trans(\vJ_i^{G_2})\trans(\vV_{it}-\vV_i\trans\vP_{it})+R_{it2}\right)(\vPi_{it}\trans\valpha_i+R_{it1})\nonumber\\
&&+2\sum_{i,t}R_{it2}\cdot(\veta-\veta_0)\trans(\vJ_i^{G_2})\trans(\vV_{it}-\vV_i\trans\vP_{it}).
\end{eqnarray}
Let $\wt\veta=\veta_0+\left(\sum_{i,t}(\vJ_i^{G_2})\trans(\vV_{it}-\vV_i\trans\vP_{it})(\vV_{it}\trans-\vP_{it}\trans\vV_i)\vJ_i^{G_2}\right)^{-1}\sum_{i,t}(\vJ_i^{G_2})\trans(\vV_{it}-\vV_i\trans\vP_{it})\epsilon_{it}$, which is actually the minimizer of the first two terms in (\ref{eqn:resid}) above. Then 
for any unit vector $\va_2\in R^{H_2}$, we have
\begin{eqnarray*}
&&\va_2\trans\vD^{G_2}(\wt\veta-\veta_0)\\
&=&T^{-1}\va_2\trans\left(\vO_2\trans\left(\begin{array}{cccc}
					\wh\vC_{11}&\vnull&\cdots&\vnull\\
					\vnull&\wh\vC_{22}&\cdots&\vnull\\					
					\vdots&\vdots&\vdots&\vdots\\
					\vnull&\vnull&\cdots&\wh\vC_{mm}
					\end{array}\right)\vO_2\right)^{-1}\cdot\\
&&	\sum_t\vO_2\trans	\left(\begin{array}{c}
						(\vV_{1t}-\vV_1\trans\vP_{1t})\epsilon_{1t}\\			
						\vdots\\
						(\vV_{mt}-\vV_m\trans\vP_{mt})\epsilon_{mt}
						\end{array}\right),
\end{eqnarray*}

Consider 
\bse
b_1&:=&T^{-1}\va_2\trans\left(\vO_2\trans\left(\begin{array}{cccc}
					\wh\vC_{11}&\vnull&\cdots&\vnull\\
					\vnull&\wh\vC_{22}&\cdots&\vnull\\					
					\vdots&\vdots&\vdots&\vdots\\
					\vnull&\vnull&\cdots&\wh\vC_{mm}
					\end{array}\right)\vO_2\right)^{-1}\cdot\\
&&	\sum_t\vO_2\trans	\left(\begin{array}{c}
						(\vV_{1t}-\vPhi_{1t})\epsilon_{1t}\\			
						\vdots\\
						(\vV_{mt}-\vPhi_{mt})\epsilon_{mt}
						\end{array}\right),\\
b_2&:=&T^{-1}\va_2\trans\left(\vO_2\trans\left(\begin{array}{cccc}
					\vC_{11}&\vnull&\cdots&\vnull\\
					\vnull&\vC_{22}&\cdots&\vnull\\					
					\vdots&\vdots&\vdots&\vdots\\
					\vnull&\vnull&\cdots&\vC_{mm}
					\end{array}\right)\vO_2\right)^{-1}\cdot\\
&&	\sum_t\vO_2\trans	\left(\begin{array}{c}
						(\vV_{1t}-\vPhi_{1t})\epsilon_{1t}\\			
						\vdots\\
						(\vV_{mt}-\vPhi_{mt})\epsilon_{mt}
						\end{array}\right).
\ese
As when showing the convergence rate, the covariance matrix of $	\left(\begin{array}{c}
						(\vV_{1t}-\vPhi_{1t})\epsilon_{1t}\\			
						\vdots\\
						(\vV_{mt}-\vPhi_{mt})\epsilon_{mt}
						\end{array}\right)$ is given by
$$\sum_{1\le t,t'\le T}\left[\begin{array}{ccc}
		\vC_{11,|t-t'|}\sigma_{11,|t-t'|}&\cdots &\vC_{1m,|t-t'|}\sigma_{1m,|t-t'|}\\
		\vdots&\vdots&\vdots\\
		\vC_{m1,|t-t'|}\sigma_{m1,|t-t'|}&\cdots &\vC_{mm,|t-t'|}\sigma_{mm,|t-t'|}\end{array}\right],$$
with eigenvalues of order $O_p(T)$ and thus $|b_2|=O_p(\sqrt{1/T})$. Using the central limit theorem under mixing conditions, for example results in \cite{bardet08}, we have
$$\sqrt{T}\nu_{2,T}^{-1/2}b_2\stackrel{d}{\rightarrow}N(0,1),$$
where 
\bse
\nu_{2,T}&=&\va_2\trans(\vO_2\trans\vC\vO_2)^{-1}\vO_2\trans\vSigma_2\vO_2(\vO_2\trans\vC\vO_2)^{-1}\va_2,\\
\vC&=&\left(\begin{array}{cccc}
					\vC_{11}&\vnull&\cdots&\vnull\\
					\vnull&\vC_{22}&\cdots&\vnull\\					
					\vdots&\vdots&\vdots&\vdots\\
					\vnull&\vnull&\cdots&\vC_{mm}
					\end{array}\right),\\
\vSigma_2&=&\frac{1}{T}\sum_{1\le t,t'\le T}\left[\begin{array}{ccc}
		\vC_{11,|t-t'|}\sigma_{11,|t-t'|}&\cdots &\vC_{1m,|t-t'|}\sigma_{1m,|t-t'|}\\
		\vdots&\vdots&\vdots\\
		\vC_{m1,|t-t'|}\sigma_{m1,|t-t'|}&\cdots &\vC_{mm,|t-t'|}\sigma_{mm,|t-t'|}\end{array}\right].
\ese
Using Lemma \ref{lem:eigenC}, $|b_1-b_2|=o_p(\sqrt{1/T})$. 
We also have that
\bse
|\va_2\trans\vD^{G_2}(\wt\veta-\veta_0)-b_1|&=&o_p(1/\sqrt{T}),
\ese
since
\bse
&&E|\va_2\trans\vD^{G_2}(\wt\veta-\veta_0)-b_1|^2\\
&=&O_p(T^{-2})\lambda_{\max}\left(E\left[\left(\begin{array}{c}
		((\vI-\vP_1)\vPhi_1-\vP_1(\vV_1-\vPhi_1))\trans\vepsilon_1\\
		\vdots\\
		((\vI-\vP_m)\vPhi_m-\vP_m(\vV_m-\vPhi_m))\trans\vepsilon_m
			\end{array}\right)^{\otimes 2}\right]\right)\\
&=&o_p(1/T).			
\ese 

Now we note that, as shown in proving convergence rate, uniformly for $\|\vtheta-\vtheta_0\|^2+\|\vbeta-\vbeta_0\|^2\le Cr_T^2$,
\be\label{eqn:claim}
&&-2\sum_{i,t}R_{it2}\epsilon_{it}\nonumber\\
&&+\sum_{i,t}R_{it2}^2+2\sum_{i,t}\left((\veta-\veta_0)(\vJ_i^{G_2})\trans(\vV_{it}-\vV_i\trans\vP_{it})+R_{it2}\right)(\vPi_{it}\trans\valpha_i+R_{it1})\nonumber\\
&&+2\sum_{i,t}R_{it2}\cdot(\veta-\veta_0)(\vJ_i^{G_2})\trans(\vV_{it}-\vV_i\trans\vP_{it})=o_p(1).
\ee
Letting 
\bse
Q(\veta)&:=&\sum_{i,t}(\veta-\veta_0)\trans(\vJ_i^{G_2})\trans(\vV_{it}-\vV_i\trans\vP_{it})(\vV_{it}\trans-\vP_{it}\trans\vV_i)\vJ_i^{G_2}(\veta-\veta_0)\\
&&-2\sum_{i,t}(\veta-\veta_0)\trans(\vJ_i^{G_2})\trans(\vV_{it}-\vV_i\trans\vP_{it})\epsilon_{it},
\ese
we have
\begin{eqnarray*}
&&\left|\sum_{i,t}(\epsilon_{it}-\vPi_{it}\trans\valpha_i-(\vV_{it}-\vV_i\trans\vP_{it})\trans(\vbeta_i-\vbeta_{0i})-R_{it1}-R_{it2})^2\right.\\
&&\left.-\sum_{i,t}(\epsilon_{it}-\vPi_{it}\trans\valpha_i-R_{it1})^2-Q(\veta)\right|=o_p(1).
\end{eqnarray*}
This implies
\begin{eqnarray*}
&&\left|\sum_{i,t}(\epsilon_{it}-\vPi_{it}\trans\valpha_i-(\vV_{it}-\vV_i\trans\vP_{it})\trans(\wb\vbeta_i-\wb\vbeta_{0i})-R_{it1}-R_{it2}(\calM_i(\valpha_i,\vbeta_i)))^2\right.\\
&&-\sum_{i,t}(\epsilon_{it}-\vPi_{it}\trans\valpha_i-(\vV_{it}-\vV_i\trans\vP_{it})\trans(\wt{\wb\vbeta}_i-\wb\vbeta_{0i})-R_{it1}-R_{it2}(\calM_i(\valpha_i,\wt\vbeta_i)))^2\\
&&\left.-(Q(\veta)-Q(\wt\veta))\right|=o_p(1).
\end{eqnarray*}
Since $Q(\veta)-Q(\wt\veta)=  \sum_{i,t}(\veta-\wt\veta)(\vJ_i^{G_2})\trans(\vV_{it}-\vV_i\trans\vP_{it})(\vV_{it}\trans-\vP_{it}\trans\vV_i)\vJ_i^{G_2}(\veta-\wt\veta)$, for any $\veta$ with $\|\vD^{G_2}(\veta-\wt\veta)\|=\delta/\sqrt{T}$ where $\delta>0$ is a small number, $Q(\veta)-Q(\wt\veta)$ is bounded away from zero. This leads to that $\sum_{i,t}(\epsilon_{it}-\vPi_{it}\trans\valpha_i-(\vV_{it}-\vV_i\trans\vP_{it})\trans(\vJ_i^{G_2})\trans(\veta-\veta_{0})-R_{it1}-R_{it2}(\calM_i (\valpha_i,\vbeta_i)))^2$ is larger than $\sum_{i,t}(\epsilon_{it}-\vPi_{it}\trans\valpha_i-(\vV_{it}-\vV_i\trans\vP_{it})\trans(\vJ_i^{G_2})\trans(\wt\veta-\veta_{0})-R_{it1}-R_{it2}(\calM_i(\valpha_i,\wt\vbeta_i)))^2$ with probability approaching one. Thus there is a local minimizer $(\wh\valpha,\wh\veta)$ of $\sum_{i,t}(\epsilon_{it}-\vPi_{it}\trans\valpha_i-(\vV_{it}-\vV_i\trans\vP_{it})\trans(\vJ_i^{G_2})\trans(\veta-\veta_{0})-R_{it1}-R_{it2}(\calM_i(\valpha_i,\vbeta_i)))^2$ with $\|\vD^{G_2}(\wh\veta-\wt\veta)\|=o_p(1/\sqrt{T})$. Thus $|\va_2\trans\vD^{G_2}(\wh\veta-\veta_0)-\va_2\trans\vD^{G_2}(\wt\veta-\veta_0)|=o_p(1/\sqrt{T})$ which proved the theorem.


Since $\wh{\wb\vbeta}-\wb\vbeta_0=\vO_2\vD^{G_2}(\wh\veta-\veta_0)$, 
$\vb_2\trans(\wh{\wb\vbeta}-\wb\vbeta_0)=\vb_2\trans\vO_2\vD^{G_2}(\wh\veta-\veta_0)$ is asymptotically normal. That is, for any unit vector $\vb_2\in R^{mp}$,
\be\label{eqn:anbeta}
\sqrt{T}\kappa_{2,T}^{-1/2}\vb_2\trans(\wh{\wb\vbeta}-\wb\vbeta_0)\stackrel{d}{\rightarrow} N(0,1),
\ee
where
$$\kappa_{2,T}:=\vb_2\trans\vO_2(\vO_2\trans\vC\vO_2)^{-1}\vO_2\trans\vSigma_2\vO_2(\vO_2\trans\vC\vO_2)^{-1}\vO_2\trans\vb_2.$$

\bigskip 

\noindent\textit{STEP 3. Proof of the convergence rate of $\wh\vtheta$ and its asymptotic normality.}

To get convergence rate of $\wh\vtheta$, like for $\wh\vbeta$, we perform a projection, which is now the projection for the nonparametric part. 
Let $\vA_{0i}:=\arg\min_{\vA}\|\vB(\vX_{it}\trans\vbeta_{0i})-g'_{0i}(\vX_{it}\trans\vbeta_{0i})\vA\wb\vX_{it}\|^2$. Obviously, we have $\vA_{0i}=E\left[g'_{0i}(\vX_{it}\trans\vbeta_{0i})\vB(\vX_{it}\trans\vbeta_{0i})\wb\vX_{it}\trans\right]\left(E\left[(g'_{0i}(\wb\vX_{it}\trans\vbeta_{0i}))^2\wb\vX_{it}\wb\vX_{it}\trans\right]\right)^{-1}$. In this part, Lemma \ref{lem:eigentheta} plays the role of assumption (C4) which was used in showing $\|\wh\vbeta-\vbeta_0\|^2=H_2/T$ previously.

Now we show $\|\wh\vtheta-\vtheta_0\|^2=O_p(H_1/T)$. The general strategy is similar to that used in showing $\|\wh\vbeta-\vbeta_0\|^2=O_p(H_2/T)$. We have
\bse
&&\sum_{i,t}(y_{it}-\vtheta_i\trans\vB(\vX_{it}\trans\vbeta_i))^2\\
&=&\sum_{i,t}(\epsilon_{it}-\vPi_{it}\trans(\vtheta_i-\vtheta_{0i})-\vV_{it}\trans(\wb\vbeta_i-\wb\vbeta_{0i})-R_{it1}-R_{it2}(\vtheta_i,\vbeta_i))^2\\
&=&\sum_{i,t}(\epsilon_{it}-(\vPi_{it}\trans-\vQ_{it}\trans\vPi_i)(\vtheta_i-\vtheta_{0i})-\vV_{it}\trans(\vgamma_i-\vgamma_{0i})-R_{it1}-R_{it2}(\vtheta_i,\vbeta_i))^2\\
&=&\sum_{i,t}(\epsilon_{it}-(\vPi_{it}\trans-\vQ_{it}\trans\vPi_i)(\vtheta_i-\vtheta_{0i})-\vV_{it}\trans(\vgamma_i-\vgamma_{0i})-R_{it1}-R_{it2}(\vtheta_i,\vgamma_i-(\vV_i\trans\vV_i)^{-1}\vV_i\trans\vtheta_i))^2,
\ese
where $\vQ_{it}\trans$ is the $t$-th row of $\vQ_i=\vV_i(\vV_i\trans\vV_i)^{-1}\vV_i\trans$ and $\vgamma_i=\wb\vbeta_i+(\vV_i\trans\vV_i)^{-1}\vV_i\trans\vPi_i\vtheta_i$, $\vgamma_{0i}=\wb\vbeta_{0i}+(\vV_i\trans\vV_i)^{-1}\vV_i\trans\vPi_i\vtheta_{0i}$.

Then
\be\label{eqn:diff2}
0
&\ge & \sum_{i,t}(\epsilon_{it}-(\vPi_{it}\trans-\vQ_{it}\trans\vPi_i)(\wh\vtheta_i-\vtheta_{0i})-\vV_{it}\trans(\wh\vgamma_i-\vgamma_{0i})-R_{it1}-R_{it2}(\wh\vtheta_i,\wh\vgamma_i-(\vV_i\trans\vV_i)^{-1}\vV_i\trans\wh\vtheta_i))^2\nonumber\\
&&-\sum_{i,t}(\epsilon_{it}-\vV_{it}\trans(\wh\vgamma_i-\vgamma_{0i})-R_{it1}-R_{it2}(\vtheta_{0i},\wh\vgamma_i-(\vV_i\trans\vV_i)^{-1}\vV_i\trans\vtheta_{0i}))^2\nonumber\\
&=&\sum_{i,t}(\wh\vxi-\vxi_0)\trans(\vJ_i^{G_1})\trans(\vPi_{it}-\vPi_i\trans\vQ_{it})(\vPi_{it}\trans-\vQ_{it}\trans\vPi_i)\vJ_i^{G_1}(\wh\vxi-\vxi_0)\nonumber\\
&&-2\sum_{i,t}(\wh\vxi-\vxi_0)\trans(\vJ_i^{G_1})\trans(\vPi_{it}-\vPi_i\trans\vQ_{it})\epsilon_{it}\nonumber\\
&&-2\sum_{i,t}(\wh R_{it2}-R_{it2})\epsilon_{it}\nonumber\\
&&+\sum_{i,t}(\wh R_{it2}^2-R_{it2}^2)\nonumber\\
&&+2\sum_{i,t}\left((\wh\vxi-\vxi_0)\trans(\vJ_i^{G_1})\trans(\vPi_{it}-\vPi_i\trans\vQ_{it})+\wh R_{it2}-R_{it2}\right)(\vV_{it}\trans(\wh\vgamma_i-\vgamma_{0i})+R_{it1})\nonumber\\
&&+2\sum_{i,t}\wh R_{it2}(\wh\vxi-\vxi_0)\trans(\vJ_i^{G_1})\trans(\vPi_{it}-\vPi_i\trans\vQ_{it}),
\ee
where we write $R_{it2}(\wh\vtheta_i,\wh\vgamma_i-(\vV_i\trans\vV_i)^{-1}\vV_i\trans\wh\vtheta_i)$ as $\wh R_{it2}$ and $R_{it2}(\vtheta_{0i},\wh\vgamma_i-(\vV_i\trans\vV_i)^{-1}\vV_i\trans\vtheta_{0i})$ as $R_{it2}$.
We have 
\bse
&&\sum_{i,t}(\wh\vxi-\vxi_0)\trans(\vJ_i^{G_1})\trans(\vPi_{it}-\vPi_i\trans\vQ_{it})(\vPi_{it}\trans-\vQ_{it}\trans\vPi_i)\vJ_i^{G_1}(\wh\vxi-\vxi_0)\\
&=&T(\wh\vxi-\vxi_0)\trans\vD^{G_1} \vO_1\trans \left(\begin{array}{cccc}
						\wh\vD_{11}&\vnull&\cdots&\vnull\\
						\vnull&\wh\vD_{22}&\cdots&\vnull\\
						\vdots&\vdots&\vdots&\vdots\\
						\vnull&\vnull&\cdots&\wh\vD_{mm}
						\end{array}\right)\vO_1\vD^{G_1} (\wh\vxi-\vxi_0)\\
&\ge & CT\|\vD^{G_1}(\wh\vxi-\vxi_0)\|^2,
\ese
where $\vO_1=\left(\begin{array}{c} \vJ_1^{G_1}\\\vdots\\\vJ_m^{G_1}\end{array}\right)(\vD^{G_1})^{-1}$ is an $mK\times H_1$ orthonormal matrix, and $\wh\vD_{ii}=\sum_{t}(\vPi_{it}-\vPi_i\trans\vQ_{it})(\vPi_{it}-\vPi_i\trans\vQ_{it})\trans/T$, and the lower bound is obtained since $\wh\vD_{ii}$ can be shown to have eigenvalues uniformly bounded from zero, similar to Lemma \ref{lem:eigenC} and using Lemma \ref{lem:eigentheta}. Furthermore, as for (\ref{eqn:second}), $\sum_{i,t}(\wh\vxi-\vxi_0)\trans(\vJ_i^{G_1})\trans(\vPi_{it}-\vPi_i\trans\vQ_{it})\epsilon_{it}=O_p(\sqrt{TH_1})$, and also the last four terms of \eqref{eqn:diff2} are $o_p(1)$, which leads to $\|\wh\vtheta-\vtheta_0\|^2=\|\vD^{G_1}(\wh\vxi-\vxi_0)\|^2=O_p(H_1/T)$.

Similarly, we can show the asymptotic normality of $\wh\vtheta$ using basically the same arguments used in showing the asymptotic normality of $\wh\vbeta$.  
Let $$\wt\vxi=\vxi_0+\left(\sum_{i,t}(\vJ_i^{G_1})\trans(\vPi_{it}-\vPi_i\trans\vQ_{it})(\vPi_{it}\trans-\vQ_{it}\trans\vPi_i)\vJ_i^{G_1}\right)^{-1}\sum_{i,t}(\vJ_i^{G_1})\trans(\vPi_{it}-\vPi_i\trans\vQ_{it})\epsilon_{it}.$$ Then 
for any unit vector $\va_1\in R^{H_1}$, we have
\begin{eqnarray*}
&&\va_1\trans\vD^{G_1}(\wt\vxi-\vxi_0)\\
&=&T^{-1}\va_1\trans\left(\vO_1\trans\left(\begin{array}{cccc}
					\wh\vD_{11}&\vnull&\cdots&\vnull\\
					\vnull&\wh\vD_{22}&\cdots&\vnull\\					
					\vdots&\vdots&\vdots&\vdots\\
					\vnull&\vnull&\cdots&\wh\vD_{mm}
					\end{array}\right)\vO_2\right)^{-1}\cdot\\
&&	\sum_t\vO_1\trans	\left(\begin{array}{c}
						(\vPi_{1t}-\vPi_1\trans\vQ_{1t})\epsilon_{1t}\\			
						\vdots\\
						(\vQ_{mt}-\vPi_m\trans\vQ_{mt})\epsilon_{mt}
						\end{array}\right),
\end{eqnarray*}
As before, it can be shown that the above is asymptotically equivalent to 
\begin{eqnarray*}
&&T^{-1}\va_1\trans\left(\vO_1\trans\left(\begin{array}{cccc}
					\vD_{11}&\vnull&\cdots&\vnull\\
					\vnull&\vD_{22}&\cdots&\vnull\\					
					\vdots&\vdots&\vdots&\vdots\\
					\vnull&\vnull&\cdots&\vD_{mm}
					\end{array}\right)\vO_2\right)^{-1}\cdot\\
&&	\sum_t\vO_1\trans	\left(\begin{array}{c}
						(\vPi_{1t}-\vPsi_{1t})\epsilon_{1t}\\			
						\vdots\\
						(\vQ_{mt}-\vPsi_{mt})\epsilon_{mt}
						\end{array}\right),
\end{eqnarray*}
where $\vD_{ii}=E[(\vB(\vX_{it}\trans\vbeta_{0i})-g'_{0i}(\vX_{it}\trans\vbeta_{0i})\vA_{0i}\wb\vX_{it})^{\otimes 2}]$ and $\vPsi_{it}=g'_{0i}(\vX_{it}\trans\vbeta_{0i})\vA_{0i}\wb\vX_{it}$. This implies that 
$$\sqrt{T}\nu_{1,T}^{-1/2}\va_1\trans\vD^{G_1}(\wt\vxi-\vxi_0)\stackrel{d}{\rightarrow} N(0,1),$$
where 
\bse
\nu_{1,T}&=&\va_1\trans(\vO_1\trans\vD\vO_1)^{-1}\vO_1\trans\vSigma_1\vO_1(\vO_1\trans\vD\vO_1)^{-1}\va_1,\\
\vD&=&\left(\begin{array}{cccc}
					\vD_{11}&\vnull&\cdots&\vnull\\
					\vnull&\vC_{22}&\cdots&\vnull\\					
					\vdots&\vdots&\vdots&\vdots\\
					\vnull&\vnull&\cdots&\vD_{mm}
					\end{array}\right),\\
\vSigma_1&=&\frac{1}{T}\sum_{1\le t,t'\le T}\left[\begin{array}{ccc}
		\vD_{11,|t-t'|}\sigma_{11,|t-t'|}&\cdots &\vD_{1m,|t-t'|}\sigma_{1m,|t-t'|}\\
		\vdots&\vdots&\vdots\\
		\vD_{m1,|t-t'|}\sigma_{m1,|t-t'|}&\cdots &\vD_{mm,|t-t'|}\sigma_{mm,|t-t'|}\end{array}\right],\\
\vD_{ii',|t-t'|}&=&E[(\vB(\vX_{it}\trans\vbeta_{0i})-g'_{0i}(\vX_{it}\trans\vbeta_{0i})\vA_{0i}\wb\vX_{it})(\vB(\vX_{i't'}\trans\vbeta_{0i'})-g'_{0i'}(\vX_{i't'}\trans\vbeta_{0i'})\vA_{0i'}\wb\vX_{i't'})\trans].		
\ese

Since $\wh{\vtheta}-\vtheta_0=\vO_1\vD^{G_1}(\wh\vxi-\vxi_0)$, 
$\vb_1\trans(\wh{\vtheta}-\vtheta_0)=\vb_1\trans\vO_1\vD^{G_1}(\wh\vxi-\vxi_0)$ is asymptotically normal. 
That is, for any unit vector $\vb_1\in R^{mp}$,
\be\label{eqn:antheta}
\sqrt{T}\kappa_{1,T}^{-1/2}\vb_1\trans(\wh{\vtheta}-\vtheta_0)\stackrel{d}{\rightarrow} N(0,1),
\ee
where
$$\kappa_{1,T}:=\vb_1\trans\vO_1(\vO_1\trans\vD\vO_1)^{-1}\vO_1\trans\vSigma_1\vO_1(\vO_1\trans\vD\vO_1)^{-1}\vO_1\trans\vb_1.$$

\subsection{Proof of Theorems \ref{thm:new1} and \ref{thm:new2} }\label{sec:conv4}
We now consider the proof of  Theorems \ref{thm:new1} and \ref{thm:new2} as special cases of (\ref{eqn:anbeta}) and (\ref{eqn:antheta}). 
Consider first Theorem \ref{thm:new2}, under the additional assumption that the true partition is used. As shown previously, the asymptotic variance of $\wh{\wb\vbeta}-\wb\vbeta_0$ is $T^{-1}\vO_2\vTheta_2\vO_2\trans$, where $\vTheta_2=(\vO_2\trans\vC\vO_2)^{-1}\vO_2\trans\vSigma_2\vO_2(\vO_2\trans\vD\vO_2)^{-1}$. From our proof, it is easy to see that eigenvalues of $\vTheta_2$ are bounded and bounded away from zero. By the definition of the $mp\times H_2$ matrix $\vO_2$, it is easy to see that its row corresponding to $\beta_{ij}$, say denoted by $\vO_{2(ij)}\trans$, has a single nonzero entry $1/\sqrt{m_{ij}}$. Let $\ve_{ij}=\sqrt{m_{ij}}\vO_{2(ij)}$, which is a unit vector, then the asymptotic variance of $\wh\beta_{ij}-\beta_{0ij}$ is $(m_{ij}T)^{-1}\ve_{ij}\trans\vTheta_2\ve_{ij}$.  

The asymptotic variance of $\wh\vtheta_i-\vtheta_{0i}$ is $T^{-1}\vJ_i^{G_1}(\vD^{G_1})^{-1}\wb\vTheta_1(\vD^{G_1})^{-1}(\vJ_i^{G_1})\trans$, where $\wb\vTheta_1=(\vO_1\trans\vD\vO_1)^{-1}\vO_1\trans\vSigma_1\vO_1(\vO_1\trans\vC\vO_1)^{-1}$ with eigenvalues bounded and bounded away from zero. By definition of $\vJ_i^{G_1}$ and $\vD^{G_1}$, it can be seen that each row of the $K\times H_1$ matrix $\vJ_i^{G_1}(\vD^{G_1})^{-1}$ has a single nonzero entry $1/\sqrt{m_i}$ and thus if we define $\vK_i=\sqrt{m_i}\vJ_i^{G_1}(\vD^{G_1})^{-1}$, it is easy to directly verify that $\vK_i\trans\vv$ is bounded and bounded away from zero and infinity for any unit vector $\vv$. Also, we have $\|\vB(x)\|\asymp K$. Thus the asymptotic variance of $\vB\trans(x)\wh\vtheta_i-\vB\trans(x)\vtheta_{0i}$ can be written as $\frac{K}{m_iT}\vb\trans(x)\vTheta_1\trans\vb(x)$, if we define $\vb(x)=\vK_i\trans\vB(x)/\|\vK_i\trans\vB(x)\|$,  and $\vTheta_1=\wb\vTheta_1\|\vK_i\trans\vB(x)\|^2/K$.

For Theorem \ref{thm:new1}, since the result is standard, and also is a special case of Theorem \ref{thm:new2}, we omit the repetition of arguments above. The  quantities $\wt\ve_{ij}$,  $\wt\vb(x)$, $\wt\vTheta_1$ and $\wt\vTheta_2$ are defined as above based on the trivial structure in which each single parameter forms its own group in the partition.

The proof of Theorem \ref{thm:new2} would be complete if we can establish consistency of homogeneity pursuit based on change point detection. That is, we need to show that the true partition can be identified with probability approaching one. 
Again for clarity the proof of this is split into three steps.

\bigskip

\noindent\textit{STEP 1. First consider the rate of $|\wh\eta_1-\eta_{01}|$.}

The proof is similar as for the rates of $\|\wh\vbeta-\vbeta_0\|$, with more complicated notations. Write $\vD^{G_2}={\rm diag}(D_1,\vD_2\}$ where $D_1=\sqrt{|G_{2,1}|}$ is the the $(1,1)$-entry of $\vD^{G_2}$, write $\vJ_i^{G_2}=(\vJ_{i1},\vJ_{i2})$ with $\vJ_{i1}$ the first column of $\vJ_i^{G_2}$.  Also write $\veta=(\eta_1,\veta_2\trans)\trans$. We have
\bse
&&\sum_{i,t}(y_{it}-\vtheta_i\trans\vB(\vX_{it}\trans\vbeta_i))^2\\
&=&\sum_{i,t}(\epsilon_{it}-\vPi_{it}\trans(\vtheta_i-\vtheta_{0i})-\vV_{it}\trans\vJ_{i2}\vD_2^{-1}\vD_2(\veta_2-\veta_{02})-\vV_{it}\trans\vJ_{i1}D_1^{-1}D_1(\eta_1-\eta_{01})  -R_{it1}-R_{it2}(\vtheta_i,\vbeta_i))^2\\
&=&\sum_{i,t}(\epsilon_{it}-\wt\vPi_{it}\trans(\vdelta_i-\vdelta_{0i})-(U_{it}-\wt\vP_{it}\trans\vU_i)D_1(\eta_1-\eta_{01})  -R_{it1}-R_{it2}(\calM_i(\vdelta_i,\eta_1))^2\\
\ese
where $\vU_{i}$ is $T$-vector with entries $U_{it}=\vV_{it}\trans\vJ_{i1}/D_1$, $\wt\vPi_{it}=(\vPi_{it}\trans,\vV_{it}\trans\vJ_2\vD_2^{-1})\trans$, $\wt\vPi_i=(\wt\vPi_{i1},\ldots,\wt\vPi_{iT})\trans$, $\wt\vP_i=\wt\vPi_i(\wt\vPi_i\trans\wt\vPi_i)^{-1}\wt\vPi_i\trans$ with rows $\wt\vP_{it}\trans=\wt\vPi_{it}(\wt\vPi_i\trans\wt\vPi_i)^{-1}\wt\vPi_i\trans$, $\vdelta=
\left(\begin{array}{c}\vtheta\\\vD_2\veta_2\end{array}\right)+(\wt\vPi_i\trans\wt\vPi_i)^{-1}\wt\vPi_i\trans \vU_iD_1\eta_1$ and $\vdelta_0=
\left(\begin{array}{c}\vtheta_0\\\vD_2\veta_{02}\end{array}\right)+(\wt\vPi_i\trans\wt\vPi_i)^{-1}\wt\vPi_i\trans \vU_iD_1\eta_{01}$. Finally, (with abuse of notation) $\calM_i(\vdelta_i,\eta_1)$ denotes the one-to-one mapping from parameterization $(\vdelta_i,\eta_1)$ to the parametrization $(\vtheta_i,\vbeta_i)$.


Then, 
\be\label{eqn:diff3}
0
&\ge & \sum_{i,t}\sum_{i,t}(\epsilon_{it}-\wt\vPi_{it}\trans(\wh\vdelta_i-\vdelta_{0i})-(U_{it}-\wt\vP_{it}\trans\vU_i)D_1(\wh\eta_1-\eta_{01})  -R_{it1}-R_{it2}(\calM_i(\wh\vdelta_i,\wh\eta_1))^2\nonumber\\
&&-\sum_{i,t}(\epsilon_{it}-\wt\vPi_{it}\trans(\wh\vdelta_i-\vdelta_{0i}) -R_{it1}-R_{it2}(\calM_i(\wh\vdelta_i,\eta_{01}))^2\nonumber\\
&=&\sum_{i,t}(D_1(\wh\eta_1-\eta_{01}))^2 (U_{it}-\wt\vP_i\trans\vU_{i})^2\nonumber\\
&&-2\sum_{i,t}D_1(\wh\eta_1-\eta_{01})(U_{it}-\wt\vP_i\trans\vU_{i})\epsilon_{it}\nonumber\\
&&-2\sum_{i,t}(\wh R_{it2}-R_{it2})\epsilon_{it}\nonumber\\
&&+\sum_{i,t}(\wh R_{it2}^2-R_{it2}^2)\nonumber\\
&&+2\sum_{i,t}\left(D_1(\wh\eta_1-\eta_{01})(U_{it}-\wt\vP_i\trans\vU_{i})+\wh R_{it2}-R_{it2}\right)(\wt\vPi_{it}\trans(\wh\vdelta_i-\vdelta_{0i})+R_{it1})\nonumber\\
&&+2\sum_{i,t}\wh R_{it2}D_1(\wh\eta_1-\eta_{01})(U_{it}-\wt\vP_i\trans\vU_{i})
\ee
with $\wh R_{it2}=R_{it2}(\calM_i(\wh\vdelta_i,\wh\eta_1))$ and $R_{it2}=R_{it2}(\calM_i(\wh\vdelta_i,\eta_{01}))$. 
The convergence rate $(\wh\eta_1-\eta_{01})^2=O_p(1/|D_1^2T|)$ is obtained by that the first term in (\ref{eqn:diff3}) is bounded below by $CTD_1^2(\wh\eta_1-\eta_{01})^2$, the second term is $O_p(\sqrt{T})|D_1(\wh\eta_1-\eta_{01})|$ while the rest are $o_p(1)$. Arguments for showing these are the same as those used in showing the rates of $\|\wh\vbeta-\vbeta_0\|$ and $\|\wh\vtheta-\vtheta_0\|$, and thus the details are omitted.

\bigskip

\noindent\textit{STEP 2. Now consider the convergence rate of $\|\vD^{G_2}(\wh\veta-\veta_0)\|_\infty$ and $\|\vD^{G_1}(\wh\vxi-\vxi_0)\|_\infty$.} 

In the study of $|\wh\eta_1-\eta_{01}|$ above, we do not make explicit that various quantities such as $U_{it}$, $\vP_i$ depends on which component of $\veta$ we are focusing on. In this section, we use subscript $(j)$, $j=1,\ldots,H_2$ to make this dependence explicit. 

To get convergence rate in infinity norm, we only need to get uniform bound for the terms in (\ref{eqn:diff3}). $\sum_{i,t}(U_{it(j)}-\wt\vP_{i(j)}\trans\vU_{i(j)})^2$ is (uniformly over different components $j$ of $\eta_i$) lower bounded by $CT$ using Lemma \ref{lem:matrixlemma} and the arguments used in Lemmas \ref{lem:eigenA} and \ref{lem:eigentheta}.

For the second term in (\ref{eqn:diff3}), using Theorem 2.19 of \cite{fanyao03}, (assuming $\epsilon_{it}$ is subGaussian)
$$\max_{1\le i\le m,1\le j\le H_2}\sum_t(U_{it(j)}-\psi_{it(j)})\epsilon_{it}=O_p(\log(Tm)\sqrt{T}).$$
Using (\ref{eqn:bd1})-(\ref{eqn:bd3}), 
$$\max_{1\le i\le m,1\le j\le H_2}\sum_t(\wt\vP_{i(j)}\trans\vU_{i(j)}-\psi_{it(j)})\epsilon_{it}=O_p(\sqrt{T}).$$

The rest terms in (\ref{eqn:diff3}) are uniformly $o_p(1)$ as shown before. These calculations combined implies and convergence rate in infinity norm. 

That $\|\vD^{G_1}(\wh\vxi-\vxi_0)\|_\infty^2=O_p(\log (Tm)/T)$ can be derived in the same way and thus omitted. 

\bigskip

\noindent\textit{STEP 3. Finally we show the consistency of change point detection}.

We use sequence $b_{(1)}\le\cdots\le b_{(n)}$ $(n=mp)$ for illustration, with estimated change points $\hat k_0=0<\hat k_1<\cdots<\hat k_{\hat H_2}=n$. The true ordered sequence of $\vbeta$ is $\beta_{0(1)}\le\cdots\le \beta_{0(n)}$ with change points $k_h$, $h=0,\ldots,H_2$. Let $\gamma_2=\min_{2\le h\le H_2}|\beta_{0(k_{h+1})}-\beta_{0(k_h)}|$ be the minimum jump size. The sup-norm convergence results established above, when specializing to the estimator in stage 1, imply that $\|\wt\vbeta-\vbeta_0\|_\infty=O_p(a_T)$ where $a_T=\sqrt{\log(Tm)/T}$. On the event $\{\|\wt\vbeta-\vbeta_0\|_\infty\le Ca_T\}$. It is easy to see that
\begin{equation}\label{eqn:dev} 
\max_{s-1<k< e}|\Delta_{s,e}(k)-\Delta_{s,e}^0(k)|\le \sqrt{n}a_T.
\end{equation}
where $\Delta_{s,e}^0(k)
=
\sqrt{\frac{(e-k)(k-s+1)}{e-s+1}}
\left|
\frac{\sum_{l=k+1}^e \beta_{0(l)}}{e-k}
-
\frac{\sum_{l=s}^{k} \beta_{0(l)}}{k-s+1}
\right|.
$

Now suppose $s-1$ and $e$ are both change points and there is at least one change point inside $(s-1,e)$. Let $\wh k=\arg\max_{s-1< k<e} \Delta_{s,e}(k)$ and $k_0=\arg\max_{s-1< k<e}\Delta_{s,e}^0(k)$. We prove consistency by way of contradiction. Suppose $\wh k$ is not one of the true change points. Then there exists some $h$ such that $\wh k\in\{k_h+1,\ldots,k_{h+1}-1\}$. From Lemma 2.2 of \cite{venkatraman92}, $\Delta_{s,e}^0(k)$ is either monotone, or decreasing and then increasing on this interval, and $\max\{\Delta_{s,e}^0(k_h),\Delta_{s,e}^0(k_{h+1})\}>\Delta_{s,e}^0(\wh k)$. Assume now $\Delta_{s,e}^0(k)$ is locally decreasing at $\wh k$ (the other case would be similar). Then we have $\Delta_{s,e}^0(k_h)>\Delta_{s,e}^0(\wh k)$ and $\Delta_{s,e}^0(k_h)$ is locally decreasing on the right side of $k_h$. Then, arguing exactly as in Lemma 2.2 of \cite{cho12}, we have $\Delta_{s,e}^0(k_h)-\Delta_{s,e}^0(k_h+1)>C\gamma_2/\sqrt{n}$. This in turn leads to $\Delta_{s,e}^0(k_0)-\Delta_{s,e}^0(\wh k)>C\gamma_2/\sqrt{n}$. Since we assumed $\sqrt{n}a_T=o(\gamma_2/\sqrt{n})$, this would lead to $\Delta_{s,e}(k_0)>\Delta_{s,e}(\wh k)$ by \eqref{eqn:dev}, a contradiction by the definition of $\wh k$. Also, in this case, it is easy to see that $\max_{s-1< k<e}\Delta_{s,e}(k)\ge  \max_{s-1< k<e}\Delta_{s,e}^0(k)-\sqrt{n}a_n\ge C\gamma_2-\sqrt{n}a_T>\delta_2$.

Now suppose still $s,e$ are both change points but there are no other change point inside $(s,e)$. In this case, using \eqref{eqn:dev}, it is easy to see that 
$\max_{s-1< k<e}\Delta_{s,e}(k)\le \sqrt{n}a_T$. 

Since we refrain from further partitioning the interval $(s,e) $ if and only if $\max_{s-1< k<e}\Delta_{s,e}(k)<\delta_2$ with $na_T<<\delta_2<<\gamma_2$, we see that the algorithm consistently identifies exactly the true change points in $\vbeta_0$.

The proof for change point detection in $\vtheta$ is the same, and the proof of Theorem \ref{thm:new2} is complete.

\subsection{Proof of Theorem \ref{thm:new3}}
For the first statement, we just need to note that $\bar\vbeta$ is the minimizer of 
$$
\min_{\va}\sum_{i=1}^m \|\vbeta_i-\va\|^2,
$$
and all $\check \vbeta_i$ are the same, thus
$$\frac{1}{mp}\sum_{i=1}^m\|\check \vbeta_i-\vbeta_i\|^2\ge \frac{1}{mp}\sum_{i=1}^m\|\vbeta_i-\bar\vbeta\|^2 \ge c.$$
 Similarly we can show the second statement. 
\section*{Appendix B. Proof of Lemmas}

\renewcommand{\thesubsection}{B.\arabic{subsection}}

\setcounter{equation}{0}

\renewcommand{\theequation}{B.\arabic{equation}}
Define matrices 
$$\vA_{ii'}=E\left[\left(\begin{array}{c}
	\vB(\vX_{it}\trans\vbeta_{i})\\
	g_{0i}'(\vX_{it}\trans\vbeta_{i})\wb\vX_{it}
	\end{array}\right)
	\left(\begin{array}{cc}
     	\vB\trans(\vX_{it}\trans\vbeta_{i}) & g_{0i}'(\vX_{it}\trans\vbeta_{i})\wb\vX_{it}\trans
	\end{array}\right)\right]
	, 1\le i,i'\le m.$$

\begin{lemma}\label{lem:eigenA}
The eigenvalues of $\vA_{ii}$ are bounded and bounded away from zero. The largest singular value (the operator norm) of $\vA_{ii'}$, $i\neq i'$, is bounded. The bounds do not depend on $(i,i')$.
\end{lemma}
\textbf{Proof of Lemma \ref{lem:eigenA}.}
By the smoothness assumption (C3), there exists $\vgamma_i\in R^{p\times K}$, with rows $\vgamma_{ij}\trans,j=1,\ldots,p$, such that
\be\label{eqn:apprdprime}
|g_{0i}'(x)E[X_{it,j}|\vX_{it}\trans\vbeta_{i}=x]-\vgamma_{ij}\trans\vB(x)|\le CK^{-2}.
\ee

We show that the operator norm of $\vgamma_i$ is bounded. If $p$ is fixed, since $\|\vgamma_{ij}\|\asymp \|\vgamma_{ij}\trans\vB(.)\|_{L_2}$ is bounded, we see the operator norm of $\vgamma_i$ is bounded since it is smaller than the operator norm. In general, we use the following more complicated arguments. Since $E[\vX_{it}\vX_{it}\trans]$ has bounded eigenvalues, so does $Var(\vX_{it})$ (the covariance matrix of $\vX_{it}$) and $(\vmu_{it})^{\otimes 2}$ where $\vmu_{it}=E[\vX_{it}]$. This implies $Var(E[\wb\vX_{it}|\vX_{it}\trans\vbeta_{i}])$ has bounded eigenvalues since $Var(\wb\vX_{it})=Var(E[\wb\vX_{it}|\vX_{it}\trans\vbeta_{i}])+E[Var(\wb\vX_{it}|\vX_{it}\trans\vbeta_{i})]$. This fact together with that $(\vmu_{it})^{\otimes 2}$ has bounded eigenvalues implies $E[(E[\wb\vX_{it}|\vX_{it}\trans\vbeta_{i}])^{\otimes 2}]$ has bounded eigenvalues. Now using (\ref{eqn:apprdprime}), $E[\vgamma_{i}\vB(\vX_{it}\trans\vbeta_{i})\vB\trans(\vX_{it}\trans\vbeta_{i})\vgamma_{i}\trans]$ has bounded eigenvalues (if $p/K^{d'}\rightarrow 0$). Since $E[\vgamma_{i}\vB(\vX_{it}\trans\vbeta_{i})\vB\trans(\vX_{it}\trans\vbeta_{i})\vgamma_{i}\trans]=\vgamma_i\vG\vgamma_i\trans$ for $\vG=E[\vB(\vX_{it}\trans\vbeta_{i})\vB\trans(\vX_{it}\trans\vbeta_{i})]$ which has eigenvalues bounded and bounded away from zero by assumption (C1). We have that the operator norm of $\vgamma_i\vG^{1/2}$ is bounded, which in turn implies the operator norm of $\vgamma_i$ is bounded. 

Then we  show that the operator norm of 
\be\label{eqn:gammai}
\left(\begin{array}{cc}
\vI &\vnull\\
-\vgamma_i &\vI\end{array}\right)
\ee
is bounded. This is easily shown by definition, since
\bse
&&\left\|\left(\begin{array}{cc}
\vI &\vnull\\
-\vgamma_i &\vI\end{array}\right)
\left(\begin{array}{c}
\vu\\
\vv\end{array}\right)\right\|\\
&=&
\left\|\left(\begin{array}{c}\vu\\ \vv-\vgamma_i\vu\end{array}\right)\right\|\\
&\le& \|\vu\|^2+2(\|\vv\|^2+\|\vgamma_i\|_{op}^2\|\vu\|^2)\le C(\|\vu\|^2+\|\vv\|^2)
\ese
Note that the inverse of (\ref{eqn:gammai}) is
$
\left(\begin{array}{cc}
\vI &\vnull\\
\vgamma_i &\vI\end{array}\right)
$
which also has bounded operator norm.

Premultiplying $\vA_{ii}$ by \eqref{eqn:gammai}
and post-multiply $\vA_{ii}$ by the transpose of \eqref{eqn:gammai},
we  get the matrix
$$E\left[\left(\begin{array}{c}
	\vB(\vX_{it}\trans\vbeta_{i})\\
	g_{0i}'(\vX_{it}\trans\vbeta_{i})\wb\vX_{it}-\vgamma_i\vB(\vX_{it}\trans\vbeta_{i})
	\end{array}\right)
	\left(\begin{array}{cc}
     	\vB\trans(\vX_{it}\trans\vbeta_{i}) & g_{0i}'(\vX_{it}\trans\vbeta_{i})\wb\vX_{it}\trans-\vB\trans(\vX_{it}\trans\vbeta_{i})\vgamma_i\trans
	\end{array}\right)\right].
$$
The operator norm for the difference between the above  and
$$E\left[\left(\begin{array}{c}
	\vB(\vX_{it}\trans\vbeta_{i})\\
	g_{0i}'(\vX_{it}\trans\vbeta_{i})(\wb\vX_{it}-E[\wb\vX_{it}|\vX_{it}\trans\vbeta_{i}])
	\end{array}\right)
	\left(\begin{array}{cc}
     	\vB\trans(\vX_{it}\trans\vbeta_{i}) & g_{0i}'(\vX_{it}\trans\vbeta_{i})(\wb\vX_{it}\trans-E[\wb\vX_{it}\trans|\vX_{it}\trans\vbeta_{i}])
	\end{array}\right)\right].
$$
is (using operator norm is bounded by the maximum row sum of absolute values of entires) $CK^{-2}(\sqrt{K}+p)=o(1)$. The displayed matrix above is block diagonal and the eigenvalues of both blocks are bounded and bounded away from zero by assumptions (C1) and (C4). This proves the first statement of the lemma.

For $\vA_{ii'}$ with $i\neq i'$, using Cauchy-Schwarz inequality, it is easy to see that for any $\vu,\vv\in R^{K+p}$, $\vu\trans\vA_{ii'}\vv\le \sqrt{\vu\trans\vA_{ii}\vu}\sqrt{\vv\trans\vA_{i'i'}\vv}$ which leads to the desired result.
\hfill $\Box$

Let $$\wh\vA_{ii'}=\frac{1}{T}\sum_{t=1}^T\left[\left(\begin{array}{c}
	\vB(\vX_{it}\trans\vbeta_{i})\\
	g_{0i}'(\vX_{it}\trans\vbeta_{i})\wb\vX_{it}
	\end{array}\right)
	\left(\begin{array}{cc}
     	\vB\trans(\vX_{it}\trans\vbeta_{i}) & g_{0i}'(\vX_{it}\trans\vbeta_{i})\wb\vX_{it}\trans
	\end{array}\right)\right]
	, 1\le i,i'\le m.$$
\begin{lemma}\label{lem:eigenAhat}
The eigenvalues of $\wh\vA_{ii}$ are bounded and bounded away from zero, and the largest singular value (the operator norm) of $\wh\vA_{ii'}$, $i\neq i'$, is bounded, with probability approaching one, uniformly over $(i,i')$ and $\vbeta$ in a neighborhood of $\vbeta_0$.
\end{lemma}
\textbf{Proof of Lemma \ref{lem:eigenAhat}.} For any $1\le k,k'\le K$ and $1\le i,i'\le m$, we have 
$$B_k(\vX_{it}\trans\vbeta_{i})B_{k'}(\vX_{i't}\trans\vbeta_{i'})\le K,$$
and 
$$E[(B_k(\vX_{it}\trans\vbeta_{i})B_{k'}(\vX_{i't}\trans\vbeta_{i'}))^2]\le KE[(B_k(\vX_{it}\trans\vbeta_{i}))^2]\le CK.$$
Thus 
$$E[(B_k(\vX_{it}\trans\vbeta_{i})B_{k'}(\vX_{i't}\trans\vbeta_{i'}))^r]\le  CK^{r-2}\cdot K,\; r=3,4,\ldots.$$
Using Theorem 2.19 of \cite{fanyao03} (setting $q=T/(C_1\log  T)$ in that theorem with large enough $C_1$), for any $\epsilon>0$,
\bse
&&P\left(\left| T^{-1}\sum_t B_k(\vX_{it}\trans\vbeta_{i})B_{k'}(\vX_{i't}\trans\vbeta_{i'})-E[B_k(\vX_{it}\trans\vbeta_{i})B_{k'}(\vX_{i't}\trans\vbeta_{i'})]\right|>\epsilon\right)\\
&\le& C(1+\log T+\mu(\epsilon))\exp\{-C\frac{T}{\log T}\mu(\epsilon)\}+CT(1+K^{C_2}/\epsilon)T^{-C_3},
\ese
where $\mu(\epsilon)=\epsilon^2/(K+K\epsilon)$, $C_2$ is some positive constant, and the constant $C_3$ can be arbitrarily large as long as one chooses $C_1$ large. Setting $\epsilon=\delta/K$,  we get 
\be
&&P(\max_{k,k',i,i'}\left| T^{-1}\sum_t B_k(\vX_{it}\trans\vbeta_{i})B_{k'}(\vX_{i't}\trans\vbeta_{i'})-E[B_k(\vX_{it}\trans\vbeta_{i})B_{k'}(\vX_{i't}\trans\vbeta_{i'})]\right|>\delta/K)\nonumber\\
&=& o(1).\label{eqn:bern1}
\ee

Similarly 
$$E\left[\left| g_{0i}'(\vX_{it}\trans\vbeta_{i})g_{0i'}'(\vX_{i't}\trans\vbeta_{i'})X_{it,j}X_{it,j'}\right|^r\right]\le  r!C^{r-2},$$
implies
\bse
&&P\left(\left|T^{-1}\sum_t g_{0i}'(\vX_{it}\trans\vbeta_{i})g_{0i'}'(\vX_{i't}\trans\vbeta_{i'})X_{it,j}X_{it,j'}-E\left[ g_{0i}'(\vX_{it}\trans\vbeta_{i})g_{0i'}'(\vX_{i't}\trans\vbeta_{i'})X_{it,j}X_{it,j'}\right]\right|>\epsilon\right)\\
&\le& C(1+\log T+\mu(\epsilon))\exp\{-C\frac{T}{\log T}\mu(\epsilon)\}+CT(1+C/\epsilon)T^{-C_3},
\ese
where $\mu(\epsilon)=\epsilon^2/(1+\epsilon)$. Setting $\epsilon=\delta/p$,  we get 
\be
&&P\left(\max_{j,j',i,i'}\left|T^{-1}\sum_t g_{0i}'(\vX_{it}\trans\vbeta_{i})g_{0i'}'(\vX_{i't}\trans\vbeta_{i'})X_{it,j}X_{it,j'}-E\left[ g_{0i}'(\vX_{it}\trans\vbeta_{i})g_{0i'}'(\vX_{i't}\trans\vbeta_{i'})X_{it,j}X_{it,j'}\right]\right|>\delta/p\right)\nonumber\\
&=& o(1).\label{eqn:bern2}
\ee

Thus $\max_{i,i'}\|\wh\vA_{ii'}-\vA_{ii'}\|_{op}=o_p(1)$ and then Lemma \ref{lem:eigenA} implies the result for any fixed $\vbeta$. 

It is easy to extend the results to obtain uniformity over $\vbeta$ in a neighborhood of $\vbeta_0$. Choosing a $T^{-a}$-covering, say $\calN_i$ of $\{\vbeta_i: \|\vbeta_i-\vbeta_{0i}\|\le b\}$ for some constant $a$ large enough. That is, for any $\vbeta_i$ there exists a $\vbeta_i'\in\calN_i$ with $\|\vbeta_i'-\vbeta_i\|<n^{-a}$. 
The size of $\calN_i$ is bounded by $\exp\{Cpa\log(T)\}$ by Lemma 2.5 of \cite{geer00}. 

To modify (\ref{eqn:bern1}) to be uniform over $\vbeta$, note that by Lipschitz continuity, it is easy to see that we have
\begin{equation}\label{eqn:lip1}
\left|T^{-1}\sum_t B_k(\vX_{it}\trans\vbeta_{i})B_{k'}(\vX_{i't}\trans\vbeta_{i'})- T^{-1}\sum_t B_k(\vX_{it}\trans\vbeta_{i}')B_{k'}(\vX_{i't}\trans\vbeta_{i'})\right|'\le T^{-a'}
\end{equation}
and 
\begin{equation}\label{eqn:lip2}
\left|E[B_k(\vX_{it}\trans\vbeta_{i})B_{k'}(\vX_{i't}\trans\vbeta_{i'})]-E[B_k(\vX_{it}\trans\vbeta_{i}')B_{k'}(\vX_{i't}\trans\vbeta_{i'}')]\right|\le T^{-a'},
\end{equation}
for some $a'>0$ (obviously we can make $a'$ arbitrarily large by setting $a$ to be large). 

Using Theorem 2.19 of \cite{fanyao03} (setting now $q=T^{1-\delta}/\log  T$), for any $\epsilon>0$,
\bse
&&P\left(\left| T^{-1}\sum_t B_k(\vX_{it}\trans\vbeta_{i})B_{k'}(\vX_{i't}\trans\vbeta_{i'})-E[B_k(\vX_{it}\trans\vbeta_{i})B_{k'}(\vX_{i't}\trans\vbeta_{i'})]\right|>\epsilon\right)\\
&\le& C(1+ T^{\delta}+\mu(\epsilon))\exp\{-CT^{1-\delta}\mu(\epsilon)\}+CT(1+K^{C_2}/\epsilon)\exp\{-C T^{-\delta}\},
\ese
where $\mu(\epsilon)=\epsilon^2/(K+K\epsilon)$. By union bound, we can still have
\bse
&&P(\max_{k,k',i,i',\vbeta_i\in\calN_i,\vbeta_{i'}\in\calN_{i'}}\left| T^{-1}\sum_t B_k(\vX_{it}\trans\vbeta_{i})B_{k'}(\vX_{i't}\trans\vbeta_{i'})-E[B_k(\vX_{it}\trans\vbeta_{i})B_{k'}(\vX_{i't}\trans\vbeta_{i'})]\right|>\delta/K)\\
&=& o(1).
\ese
 The uniformly of $\vbeta_i\in\calN_i$ imply the uniformity of $\vbeta_i$ in a neighborhood of $\vbeta_{0i}$ by (\ref{eqn:lip1}) and (\ref{eqn:lip2}).

Similarly we can modify (\ref{eqn:bern2}) to be uniform over $\vbeta$ if $p^3(\log T)^2\log (pm)/T\rightarrow 0$, which finishes the proof.
\hfill $\Box$

\begin{lemma}\label{lem:eigenAtilde}
Eigenvalues of $\wt\vA_{ii'}$ are bounded and bounded away from zero, with probability approaching one, uniformly over $(i,i')$ and $\vbeta$.
\end{lemma}
\textbf{Proof of Lemma \ref{lem:eigenAtilde}.} 
First, by Lemma \ref{lem:eigenAhat}, the eigenvalues of $(1/T)\sum_t \vB(\vX_{it}\trans\vbeta_i)\vB\trans(\vX_{it}\trans\vbeta_i)$ are bounded and bounded away from zero. 

Since $\|\vbeta_i-\vbeta_{0i}\|=O(r_T)$ and $\|\vbeta_i^*-\vbeta_{0i}\|=O(r_T)$, 
\begin{eqnarray*}
&&\frac{1}{T}\sum_{t=1}^T\left(\vtheta_{0i}\trans\vB(\vX_{it}\trans\vbeta_i^*)\right)^2X_{it,j}X_{it,j}-(g_{0i}(\vX_{it}\trans\vbeta_{i}))^2X_{it,j}X_{it,j'}\\
&=&\frac{1}{T}\sum_{t=1}^T\left(\left(\vtheta_{0i}\trans\vB(\vX_{it}\trans\vbeta_i^*)\right)^2-\left(\vtheta_{0i}\trans\vB(\vX_{it}\trans\vbeta_{i})\right)^2\right)X_{it,j}X_{it,j'}+\left(\left(\vtheta_{0i}\trans\vB(\vX_{it}\trans\vbeta_i^*)\right)^2-(g_{0i}(\vX_{it}\trans\vbeta_{i}))^2\right)X_{it,j}X_{it,j'}\\
&=&O_p(r_T\sqrt{p}+K^{-2}),
\end{eqnarray*}
\begin{eqnarray*}
&&\frac{1}{T}\sum_{t=1}^T\vtheta_{0i}\trans\vB(\vX_{it}\trans\vbeta_{i}^*) X_{it,j}B_k(\vX_{it}\trans\vbeta_i)-g_{0i}(\vX_{it}\trans\vbeta_{i}) X_{it,j}B_k(\vX_{it}\trans\vbeta_{i})\\
&=&\frac{1}{T}\sum_{t=1}^T (\vtheta_{0i}\trans\vB(\vX_{it}\trans\vbeta_{i})-g_{0i}(\vX_{it}\trans\vbeta_{i})) X_{it,j}B_k(\vX_{it}\trans\vbeta_{i})+ (\vtheta_{0i}\trans\vB(\vX_{it}\trans\vbeta_{i}^*)-\vtheta_{0i}\trans\vB(\vX_{it}\trans\vbeta_{i})) X_{it,j}B_k(\vX_{it}\trans\vbeta_{i})\\
&=&O_p((K^{-2}+r_T\sqrt{p})/K).
\end{eqnarray*}
Thus if $r_T\sqrt{p^3}=o(1)$ and $p=o(K^d)$, we have $\|\wt\vA_{ii}-\wh\vA_{ii}\|_{op}=o_p(1)$ which proves the lemma.\hfill $\Box$

\begin{lemma}\label{lem:eigenAtilde2}
Eigenvalues of 
$$\sum_{t,t'}\left[\begin{array}{ccc}
		\vA_{11,|t-t'|}\sigma_{11,|t-t'|}&\cdots &\vA_{1m,|t-t'|}\sigma_{1m,|t-t'|}\\
		\vdots&\vdots&\vdots\\
		\vA_{m1,|t-t'|}\sigma_{m1,|t-t'|}&\cdots &\vA_{1m,|t-t'|}\sigma_{1m,|t-t'|}\end{array}\right]$$
		are bounded by $CT$ for some constant $C$.
\end{lemma}
\textbf{Proof of Lemma \ref{lem:eigenAtilde2}.} By Lemma \ref{lem:eigenA}, and similar to the proof of Lemma \ref{lem:eigenAtilde}, $\max_i\|\vA_{ii,0}\|_{op}$ is bounded. Using Cauchy-Schwarz inequality, it is easy to show $\|\vA_{ii',|t-t'|}\|_{op}\le (\|\vA_{ii,0}\|_{op}+\|\vA_{i'i',0}\|_{op})/2$ and thus $\max_{i,i',t,t'}\|\vA_{ii',|t-t'|}\|_{op}$ is also bounded.

Let $\vv=(\vv_1\trans,\ldots,\vv_m\trans)\trans\in R^{m(K+p)}$.  We have
\begin{eqnarray*}
&&\sum_{t,t'}\vv\trans 
\left(\begin{array}{ccc}
		\vA_{11,|t-t'|}\sigma_{11,|t-t'|}&\cdots &\vA_{1m,|t-t'|}\sigma_{1m,|t-t'|}\\
		\vdots&\vdots&\vdots\\
		\vA_{m1,|t-t'|}\sigma_{m1,|t-t'|}&\cdots &\vA_{mm,|t-t'|}\sigma_{mm,|t-t'|}\end{array}\right)\vv\\
&=&\sum_{t,t'}\sum_{i,i'}\sigma_{ii',|t-t'|}\vv_i\trans\vA_{ii',|t-t'|}\vv_{i'}\\
&=&T\sum_{l=1}^{T-1}\sum_{i,i'}(1-l/T)\sigma_{ii',l}\vv_i\trans\vA_{ii',l}\vv_{i'}\\
&\le&T\sum_{l=1}^{T-1}\sum_{i,i'}(1-l/T)|\sigma_{ii',l}|\cdot\|\vv_i\|\cdot\|\vv_{i'}\|\cdot\lambda_{\max}(\vA_{ii',l})\\
&\le&CT\sum_{i,i'}\tau_{ii'}\cdot\|\vv_i\|\cdot\|\vv_{i'}\|\\
&\le&CT \lambda_{\max}(\{\tau_{ii'}\}_{i,i'=1}^m)\le CT,
\end{eqnarray*}
where in the last step we used that $\lambda_{\max}(\{\tau_{ii'}\}_{i,i'=1}^m)$ is  bounded, by assumption (C2) and the Gershgorin circle theorem. \hfill $\Box$

\begin{lemma}\label{lem:eigenC} 
$\max_{1\le i\le m}\|\wh\vC_{ii}-\vC_{ii}\|_{op}=o_p(1).$
\end{lemma}
\textbf{Proof of Lemma \ref{lem:eigenC}.} Let $\wt\vC_{ii}=T^{-1}\sum_{t=1}^T(g_{0i}'(\vX_{it}\trans\vbeta_{0i}))^2(\wb\vX_{it}-E[\wb\vX_{it}|\vX_{it}\vbeta_{0i}])^{\otimes 2}$. Also let $\vPhi_{it}=g_{0i}'(\vX_{it}\trans\vbeta_{0i})E[\wb\vX_{it}|\vX_{it}\trans\vbeta_{0i}]$, $\vPhi_i=(\vPhi_{i1},\ldots,\vPhi_{iT})\trans$. We have 
\begin{eqnarray*}
&&\wh\vC_{ii}-\wt\vC_{ii}\\
&=&\frac{1}{T}\sum_t(\vV_{it}-\vV_i\trans\vP_{it})(\vV_{it}\trans-\vP_{it}\trans\vV_i)-\frac{1}{T}\sum_t(\vV_{it}-\vPhi_{it})(\vV_{it}-\vPhi_{it})\trans\\
&=&\frac{1}{T}\vV_i\trans(\vI-\vP_i)\vV_i-\frac{1}{T}(\vV_i-\vPhi_i)\trans(\vV_i-\vPhi_i).
\end{eqnarray*}
Writing $\vV_i=(\vV_i-\vPhi_i)+\vPhi_i$, the above is equal to
\begin{eqnarray}\label{eqn:diff} 
&&\frac{1}{T}\left(\vPhi_i\trans(\vI-\vP_i)\vPhi_i+(\vV_i-\vPhi_i)\trans(\vI-\vP_i)\vPhi_i+\vPhi_i(\vI-\vP_i)(\vV_i-\vPhi_i)-(\vV_i-\vPhi_i)\trans\vP_i(\vV_i-\vPhi_i)\right).\nonumber\\
\end{eqnarray}

Since $g_{0i}'(\vX_{it}\trans\vbeta_{0i})E[\wb\vX_{it}|\vX_{it}\trans\vbeta_{0i}]$ is a $d'$-smooth function of $\vX_{it}\trans\vbeta_{0i}$, we have 
\begin{equation}
\label{eqn:bd1}
\|(\vI-\vP_i)\vPhi_i\|\le C\sqrt{T}K^{-2}.
\end{equation}
We also have trivially 
\begin{equation}
\label{eqn:bd2}
\max_i\|\vV_i-\vPhi_i\|=O_p(\sqrt{Tp}).
\end{equation}

Now consider $\|\vP_i(\vV_i-\vPhi_i)\|$. 
We have
\begin{eqnarray*}
&&\|\vP_i(\vV_i-\vPhi_i)\|\\
&=&\|\vPi_i(\vPi_i\trans\vPi_i)^{-1}\vPi_i\trans(\vV_i-\vPhi_i)\|\\
&\le&\|\vPi_i(\vPi_i\trans\vPi_i)^{-1}\|_{op}\|\vPi_i\trans(\vV_i-\vPhi_i)\|,
\end{eqnarray*}
and $\max_i\|\vPi_i(\vPi_i\trans\vPi_i)^{-1}\|_{op}^2=\max_i\|(\vPi_i\trans\vPi_i)^{-1}\|_{op}=O_p(1/T)$ as proved in Lemma \ref{lem:eigenAhat}. 
For the term $\|\vPi_i\trans(\vV_i-\vPhi_i)\|$, we can deal with it similar to (\ref{eqn:bern1}). First note that $B_k(\vX_{it}\trans\vbeta_{0i})(X_{it,j}-\phi_{it,j})$ has mean zero ($\phi_{it,j}$ is the $j$-th component of $\vPhi_{it}$, $j=1,\ldots,p$). Since $X_{it,j}$ is bounded, we have
\bse
&&B_k(\vX_{it}\trans\vbeta_{0i})(X_{it,j}-\phi_{it,j})\le C\sqrt{K}\\
&&E[(B_k(\vX_{it}\trans\vbeta_{0i})(X_{it,j}-\phi_{it,j}))^2]\le C,
\ese
and applying Theorem 2.19 of \cite{fanyao03}, 
\bse
&&P(B_k(\vX_{it}\trans\vbeta_{0i})(X_{it,j}-\phi_{it,j})>T\epsilon)\\
&\le &C(1+\log T+\mu(\epsilon))\exp\{-C\frac{T}{\log T}\mu(\epsilon)\}+CT(1+1/\epsilon)T^{-C_2}
\ese
where $\mu(\epsilon)=\frac{\epsilon^2}{\sqrt{K}\epsilon+1}$. Setting $\epsilon=\sqrt{\log T/T}$, and taking union bound over $(i,j,k)$, we can obtain 
\begin{equation*}
\max_i\|\vPi_i\trans(\vV_i-\vPhi_i)\|=O_p(\sqrt{TKp\log T}).
\end{equation*}
Thus  
\begin{equation}
\label{eqn:bd3}
\|\vP_i(\vV_i-\vPhi_i)\|=o_p(\sqrt{T}) 
\quad
\mbox{if } 
Kp\log T/T\rightarrow 0. 
\end{equation}

Then using (\ref{eqn:bd1})-(\ref{eqn:bd3}), (\ref{eqn:diff})  is $o_p(1)$. Finally, using the same arguments as in  the proof of Lemma \ref{lem:eigenAhat},  we get $\max_i\|\wt\vC_{ii}-\vC_{ii}\|_{op}=o_p(1)$.
 \hfill $\Box$

\begin{lemma}\label{lem:eigentheta}
Eigenvalues of $E\left[(\vB(\vX_{it}\trans\vbeta_{0i})-g'_{0i}(\vX_{it}\trans\vbeta_{0i})\vA_{0i}\wb\vX_{it})(\vB(\vX_{it}\trans\vbeta_{0i})-g'_{0i}(\vX_{it}\trans\vbeta_{0i})\vA_{0i}\wb\vX_{it})\trans\right]$ are bounded and bounded away from zero, uniformly over $i$.
\end{lemma}
The proof is based on the following elementary lemma.
\begin{lemma}\label{lem:matrixlemma}
Suppose a positive definite matrix $\left(\begin{array}{cc} \vA&\vB\\\vB\trans&\vC\end{array}\right)$ has all eigenvalues inside the interval $[c,C]$ for some $0<c<C<\infty$. 
Then all eigenvalues of $\vC-\vB\trans\vA^{-1}\vB$ are also inside the interval $[c,C]$.
\end{lemma}
\textbf{Proof of Lemma \ref{lem:matrixlemma}.} Obviously eigenvalues of $\vC-\vB\trans\vA^{-1}\vB$ are no larger than that of $\vC$, which is in turn bounded by $C$. Next, we have the identity
$$
\left(\begin{array}{cc} \vI&\vnull\\-\vB\trans\vA^{-1}&\vI\end{array}\right)
\left(\begin{array}{cc} \vA&\vB\\\vB\trans&\vC\end{array}\right)
\left(\begin{array}{cc} \vI& -\vA^{-1}\vB\\\vnull&\vI\end{array}\right)=
\left(\begin{array}{cc} \vA&\vnull\\\vnull\trans&\vC-\vB\trans\vA^{-1}\vB\end{array}\right)
$$
Thus for any vector $\vb$ with dimension same as that of $\vC$, we have
\bse
&&\vb\trans(\vC-\vB\trans\vA^{-1}\vB)\vb\\
&=&(\vnull\trans,\vb\trans)\left(\begin{array}{cc} \vA&\vnull\\\vnull\trans&\vC-\vB\trans\vA^{-1}\vB\end{array}\right)\left(\begin{array}{c} \vnull\\\vb\end{array}\right)\\
&=&(\vnull\trans,\vb\trans)
\left(\begin{array}{cc} \vI&\vnull\\-\vB\trans\vA^{-1}&\vI\end{array}\right)
\left(\begin{array}{cc} \vA&\vB\\\vB\trans&\vC\end{array}\right)
\left(\begin{array}{cc} \vI& -\vA^{-1}\vB\\\vnull&\vI\end{array}\right)
\left(\begin{array}{c} \vnull\\\vb\end{array}\right)\\
&=&(-\vb\trans\vB\trans\vA^{-1},\vb\trans)
\left(\begin{array}{cc} \vA&\vB\\\vB\trans&\vC\end{array}\right)
\left(\begin{array}{c} -\vA^{-1}\vB\vb\\\vb\end{array}\right)\\
&\ge&\|\vb\|^2 c,
\ese
which completes the proof. 
\hfill $\Box$

\noindent\textbf{Proof of Lemma \ref{lem:eigentheta}.} Since we have 
\bse
&&E\left[(\vB(\vX_{it}\trans\vbeta_{0i})-g'_{0i}(\vX_{it}\trans\vbeta_{0i})\vA_{0i}\wb\vX_{it})(\vB(\vX_{it}\trans\vbeta_{0i})-g'_{0i}(\vX_{it}\trans\vbeta_{0i})\vA_{0i}\wb\vX_{it})\trans\right]\\
&=&E[\left(\vB(\vX_{it}\trans\vbeta_{0i})\vB\trans(\vX_{it}\trans\vbeta_{0i})\right]\\
&&-
E\left[g'_{0i}(\vX_{it}\trans\vbeta_{0i})\vB(\vX_{it}\trans\vbeta_{0i})\wb\vX_{it}\trans\right]
\left(E\left[(g'_{0i}(\wb\vX_{it}\trans\vbeta_{0i}))^2\wb\vX_{it}\wb\vX_{it}\trans\right]\right)^{-1}
E\left[g'_{0i}(\vX_{it}\trans\vbeta_{0i})\wb\vX_{it}\vB\trans(\vX_{it}\trans\vbeta_{0i})\right],
\ese
the lemma  follows from Lemma \ref{lem:eigenA} and Lemma \ref{lem:matrixlemma}.\hfill $\Box$

\bibliographystyle{dcu}
\bibliography{papers,books}

\end{document}